\ifx\preambleloaded\undefined

\documentclass[11pt,a4paper]{article} %
\usepackage{fullpage} %
\usepackage{mathtools} %
\usepackage{amssymb} %
\usepackage{latexsym} %
\usepackage{amsfonts} %
\usepackage{mathrsfs} %
\usepackage{amsthm} %
\theoremstyle{plain} %

\newtheorem{thm}{Theorem}[section] %
\newtheorem{cor}[thm]{Corollary} %
\newtheorem{lem}[thm]{Lemma} %
\newtheorem{prop}[thm]{Proposition} %
\newtheorem{con}[thm]{Conjecture} %
\theoremstyle{definition} %
\newtheorem{defn}[thm]{Definition} %
\newtheorem{example}[thm]{Example} %
\newtheorem{problem}[thm]{Problem} %
\newtheorem*{notation}{Notation} %
\theoremstyle{remark} %
\newtheorem*{rem}{Remark} %

\theoremstyle{definition} %
\theoremstyle{remark} %

\makeatletter%
\renewenvironment{proof}[1][\proofname]%
{\par\pushQED{\qed}\normalfont\topsep6\p@\@plus6\p@\relax\trivlist\item[\hskip\labelsep\bfseries#1\@addpunct{.}]\ignorespaces}%
{\popQED\endtrivlist\@endpefalse}%
\makeatother %
\usepackage{graphicx} %
\usepackage{epsfig} %
\usepackage[noprefix,refpage]{nomencl} %
\makenomenclature
\usepackage[numbers]{natbib} 
\usepackage[boxed]{algorithm} %
\usepackage{algorithmic} %
\usepackage{color} %
\definecolor{darkblue}{rgb}{0.0,0.0,0.3} %
\usepackage[bookmarks,backref]{hyperref} %
\hypersetup{colorlinks,breaklinks,
            linkcolor=darkblue,urlcolor=darkblue,
            anchorcolor=darkblue,citecolor=darkblue} %

\newcommand\doilink[1]{\href{http://dx.doi.org/#1}{#1}} %
\newcommand\doi[1]{doi:\doilink{#1}} %

\providecommand*\url[1]{\href{#1}{#1}} %

\renewcommand*\url[1]{\href{#1}{\texttt{#1}}} %

\newcommand{\zz}{\mathbb{Z}} %
\newcommand{\nn}{\mathbb{N}} %
\newcommand{\card}[1]{\ensuremath{\left|#1\right|}} %
\newcommand{\tup}[2]{(#1,\dotsc,#2)}%
\newcommand{\set}[2]{\{#1,\dotsc,#2\}}%
\newcommand{\ints}[2]{[#1,#2]}%
\newcommand{\seq}[2]{#1,\dotsc,#2}%

\newcommand{\inj}{\to} %
\newcommand{\sur}{\overset{\mathtt{onto}}{\longrightarrow}} %
\newcommand{\bij}{\to} %

\makeatletter
\newcommand*{\@old@slash}{}\let\@old@slash\slash
\def\slash{\relax\ifmmode\delimiter"502F30E\mathopen{}\else\@old@slash\fi} %
\makeatother
\def\backslash{\delimiter"526E30F\mathopen{}} %

\newcommand{\iso}[1]{I(#1)} %
\newcommand{\rep}[1]{#1^\ast} %
\newcommand{\lab}{\mcg} %
\newcommand{\labc}{\mcg^c} %
\newcommand{\unlab}{\mcg/{\cong}} %
\DeclareMathOperator{\mon}{mon} %
\newcommand{\rhs}{right side } %
\newcommand{\mcc}{\mathcal{C}} %
\newcommand{\mcf}{\mathcal{F}} %
\newcommand{\mcg}{\mathcal{G}} %
\newcommand{\mch}{\mathcal{H}} %
\newcommand{\mcp}{\mathcal{P}} %
\newcommand{\mcq}{\mathcal{Q}} %
\newcommand{\mcr}{\mathcal{R}} %
\newcommand{\mcs}{\mathcal{S}} %

\newif\ifenglish
\englishtrue %
\newif\ifportu
\ifportu
\usepackage[english,brazil]{babel} %
\else
\usepackage[brazil,english]{babel} %
\fi
\usepackage[utf8]{inputenc} %
\newif\ifcurrent
\currenttrue %
\newif\ifold %

\newcommand{\mc}[1]{\leavevmode\marginpar{%
\raggedright\footnotesize
{\bf\itshape\hrule\smallskip#1\smallskip\hrule}}} %
\newcommand{\removemc}{\renewcommand{\mc}[1]{}} %

\DeclareMathOperator{\parts}{\Lambda} %
\DeclareMathOperator{\vsubop}{ind} %
\newcommand{\vsub}[2]{\vsubop(#1,#2)} %
\DeclareMathOperator{\esubop}{esub} %
\newcommand{\esub}[2]{\esubop(#1,#2)} %
\DeclareMathOperator{\subop}{sub} %
\newcommand{\sub}[2]{\subop(#1,#2)} %
\DeclareMathOperator{\covop}{cov} %
\newcommand{\cov}[2]{\covop(#1\rightarrow #2)} %
\newcommand{\vcov}[2]{\covop(#1\xrightarrow{v}#2)} %

\DeclareMathOperator{\deckop}{deck} %
\newcommand{\deck}[1]{\deckop{(#1)}} %
\newcommand{\eplus}[1]{\{#1^{+e}\}} %
\newcommand{\eminus}[1]{#1^{-e}} %
\newcommand{\down}[2]{#1(#2)} %
\newcommand{\vpo}{\leq_v} %
\newcommand{\ispu}[1]{\mcp(#1)} %
\newcommand{\isp}[1]{\overline{\mcp}(#1)} %
\newcommand{\ispug}{\mcp} %
\newcommand{\wv}{w_v} %
\newcommand{\latpo}{\vdash_c} %
\newcommand{\modpo}{\models_c} %
\newcommand{\latl}[1]{\Pi^c_{#1}} %
\newcommand{\latu}[1]{\Omega(#1)} %
\newcommand{\lat}[1]{\overline{\Omega}(#1)} %
\newcommand{\latug}{\Omega} %
\newcommand{\wl}{w_\pi} %
\newcommand{\epo}{\leq_e} %
\newcommand{\espu}[1]{\mcq(#1)} %
\newcommand{\esp}[1]{\overline{\mcq}(#1)} %
\newcommand{\espug}{\mcq} %
\newcommand{\we}{w_e} %
\usepackage[inline]{enumitem} %
\newlist{pcases}{enumerate}{1} %
\setlist[pcases]{ %
  label=\underline{Case~\arabic*:}\protect\thiscase.~, %
  ref=\arabic*, %
  align=left, %
  labelsep=0pt, %
  leftmargin=0pt, %
  labelwidth=0pt, %
  parsep=0pt %
} %
\newcommand{\pcase}[1][]{%
  \if\relax\detokenize{#1}\relax %
    \def\thiscase{} %
  \else %
    \def\thiscase{~#1} %
  \fi %
  \item %
} %
\makeatletter

\newlist{pclaims}{enumerate}{1} %
\setlist[pclaims]{ %
  label={\em Claim~\arabic*.}\protect\thisclaim~, %
  ref=\arabic*, %
  align=left, %
  labelsep=0pt, %
  leftmargin=0pt, %
  labelwidth=0pt, %
  listparindent=\parindent %
} %
\newcommand{\pclaim}[1][]{ %
  \if\relax\detokenize{#1}\relax
    \def\thisclaim{} %
  \else
    \def\thisclaim{~#1} %
  \fi
  \item
} %
\makeatletter
\newcommand{\ppf}{\noindent{\em Proof.} } %
\newcommand{\pqed}{$\square$} %

\def\preambleloaded{} %

\fi

\title{Subgraph posets and graph reconstruction}

\author{%
  Bhalchandra D. Thatte \\
  Departamento de Matem\'atica, \\
  Universade Federal de Minas Gerais, Brasil \\
  \texttt{thatte@ufmg.br} }

\date{\today}
\begin{document}
\removemc
\maketitle
\begin{abstract} 
  We consider only finite simple undirected graphs in this paper. Let
  $G$ be an arbitrary graph. Let $\ispu{G}$ be the set consisting of
  $K_1$ and the distinct unlabelled nonempty induced subgraphs of
  $G$. The {\em abstract induced subgraph poset} of $G$ is the
  isomorphism class of the weighted poset
  $(\ispu{G},{\vpo,} \wv\colon\ispu{G}\times\ispu{G} \to \nn)$, where
  for $G_i,G_j \in \ispu{G}$ we define $G_i \vpo G_j$ if $G_i$ is an
  induced subgraph of $G_j$, and $\wv(G_i,G_j)$ is the number of induced
  subgraphs of $G_j$ that are isomorphic to $G_i$. We write $\isp{G}$
  for the isomorphism class of the weighted poset defined above. In an
  earlier paper, we showed that several invariants of $G$ can be
  computed from the abstract poset $\isp{G}$, i.e., the deck of $G$ is
  not required. In this paper, we study reconstruction questions on two
  analogously defined posets: the abstract weighted lattice $\lat{G}$ of
  distinct unlabelled connected partitions of $G$, which we call the
  {\em abstract bond lattice} of $G$, and the abstract weighted poset
  $\esp{G}$ of distinct unlabelled edge-subgraphs of $G$, which we call
  the {\em abstract edge-subgraph poset} of $G$.

  We show that $\lat{G}$ can be constructed from $\isp{G}$, and that
  $\isp{G}$ can be constructed from $\lat{G}$ if $G$ is not a star or a
  disjoint union of edges and has no isolated vertices. The first
  construction implies that if a graph invariant can be computed from
  $\lat{G}$, then it can also be computed from $\isp{G}$. An examples of
  such an invariant is the chromatic symmetric function. Since every
  tree $T$ on 2 or more vertices can be reconstructed up to isomorphism
  from $\isp{T}$, the second construction implies that every tree $T$ on
  2 or more vertices that is not a star can be reconstructed up to
  isomorphism from $\lat{T}$. We also give simple proofs that the
  chromatic symmetric function $X_{G}(x)$ and the symmetric Tutte
  polynomial $X_{G}(x;t)$ of $G$ can be computed from $\isp{G}$. The
  main tools that we use to prove these results are a generalisation to
  abstract induced subgraph posets of a lemma of Kocay in graph
  reconstruction theory and other related subgraph counting identities.

  Stanley has asked if every tree $T$ is determined up to isomorphism by
  its chromatic symmetric function $X_T(x)$. Analogously, Noble and
  Welsh have asked if every tree $T$ is determined up to isomorphism by
  its symmetric Tutte polynomial $X_T(x;t)$. We show that the two
  questions are equivalent by showing that, for every tree $T$,
  $X_T(x;t)$ is determined by $X_T(x)$.

  In Section~\ref{sec-esp}, we consider the problem of reconstructing an
  arbitrary graph $G$ up to isomorphism from its abstract edge-subgraph
  poset $\esp{G}$, which we call the $Q$-reconstruction problem, and
  study its relation to the edge reconstruction conjecture of Harary. We
  present an infinite family of graphs that are not $Q$-reconstructible,
  and show that the edge reconstruction conjecture is true if and only
  if the graphs in the family are the only graphs that are not
  $Q$-reconstructible.

  Let $\lab$ be the set of all graphs, and let $\unlab$ be the set of
  all unlabelled graphs (isomorphism classes). Let $\hom(G,H)$ denote
  the number of homomorphisms from $G$ to $H$. Let
  $f\colon \unlab \to \unlab $ be a bijection such that for all
  $G, H \in \unlab$, we have $\hom(G,H) = \hom(f(G),f(H))$. We
  conjecture that $f(G)=G$ for all $G\in \unlab$. Our conjecture is
  motivated by Lov\'asz's homomorphism cancellation laws. We prove that
  the conjecture stated above is weaker than the edge reconstruction
  conjecture.

\end{abstract}

\tableofcontents

\section{Introduction}
\label{sec:intro}
We consider only finite undirected simple graphs in this paper.
A well-known conjecture of \citet{ulam1960} and \citet{kelly1942}, known
as the {\em vertex reconstruction conjecture} or {\em Ulam's
  conjecture}, states that every graph on 3 or more vertices is
determined up to isomorphism by its {\em deck} (the collection or {\em
  multiset} of its unlabelled vertex-deleted subgraphs).
An analogous conjecture, known as the {\em edge reconstruction
  conjecture}, was proposed by \citet{harary1964}. It states that every
graph with at least 4 edges is determined up to isomorphism by its {\em
  edge-deck} (the collection of its unlabelled edge-deleted subgraphs).
These are some of the foremost unsolved problems in graph theory. We
refer the reader to a survey of these conjectures by \citet{bondy1991}.

\subsection{Notation}
\label{sec:notation}
\subsubsection{Miscellaneous notation}
We denote the set of integers, the set of positive integers, and the set
of natural numbers (including 0) by $\zz$, $\zz^+$, and $\nn$,
respectively.  \nomenclature[4z]{$\zz$}{set of integers} %
\nomenclature[4z]{$\zz^+$}{set of positive integers} %
\nomenclature[4n]{$\nn$}{set of natural numbers (including 0)} %
We denote the family of $k$-element subsets of a set $S$ by
$\binom{S}{k}$, %
\nomenclature[4sk]{$\binom{S}{k}$}{family of $k$-element subsets of
  $S$} %
the set of $k$-element tuples (or the set of sequences of length $k$)
from $S$ by $S^k$,%
\nomenclature[4sk]{$S^k$}{set of $k$-element tuples of elements in
  $S$} %
and the powerset of $S$ by $2^S$. %
\nomenclature[22]{$2^S$}{powerset of $S$}%
When the range of an index is unspecified, e.g., as in $\sum_ia_i$, or
$\bigcup_iH_i$ or in expressions such as ``... for all $i$'', we
understand that full range of the index over which the objects in the
context are defined is implied. We use this convention especially
outside displayed mathematics.

\subsubsection{Graphs}
\label{sec:graphs}
We denote the set of all graphs by $\lab$ %
\nomenclature[4g1]{$\lab $}{set of all labelled graphs} %
and the set of connected graphs by $\labc$. %
\nomenclature[4g1c]{$\labc $}{set of all labelled connected graphs} %
Throughout this paper, we take $G$ and $H$ to be arbitrary graphs. We
denote the vertex set of $G$ by $V(G)$, %
\nomenclature[4vg]{$V(G)$}{vertex set of $G$} its edge set by $E(G)$, %
\nomenclature[4eg]{$E(G)$}{edge set of $G$} %
number of vertices in $G$ by $\nu(G)$, %
\nomenclature[5n]{$\nu(G)$}{number of vertices of $G$} %
the number of edges in $G$ by $\epsilon(G)$, %
\nomenclature[5e]{$\epsilon(G)$}{number of edges of $G$} %
and the number of components of $G$ by $c(G)$. %
\nomenclature[4c3g]{$c(G)$}{number of components of $G$} %
An {\em empty graph} is a graph with empty edge set. A \emph{null graph}
$\Phi$ is a graph with no vertices. %
\nomenclature[5z1]{$\Phi$}{null graph} %
For $X \subseteq V(G)$, we denote the subgraph of $G$ induced by $X$ by
$G[X]$, the subgraph of $G$ induced by $V(G)\setminus X$ by $G-X$, or
simply $G-u$ if $X=\{u\}$. %
\nomenclature[4g2x]{$G[X], X\subseteq V(G)$}{subgraph of $G$ induced by
  $X$} %
\nomenclature[4g2x]{$G-X, X\subseteq V(G)$}{subgraph of $G$ induced by
  $V(G)\setminus X$} %
\nomenclature[4g2u]{$G-u, u\in V(G)$}{subgraph of $G$ induced by
  $V(G)\setminus \{u\}$} %
For $E \subseteq E(G)$, we denote the subgraph of $G$ induced by $E$ by
$G[E]$, the spanning subgraph of $G$ with edge set $E$ by $G_E$, and the
spanning subgraph of $G$ with edge set $E(G)\backslash E$ by $G-E$ (or
just $G-e$ if $E=\{e\}$. %
\nomenclature[4g2e]{$G[E], E\subseteq E(G)$}{subgraph of $G$ induced by
  $E$} %
\nomenclature[4g2e]{$G-E, E\subseteq E(G)$}{spanning subgraph of $G$ with
  edge set $E(G)\setminus E$} %
\nomenclature[4g2e]{$G_E, E\subseteq E(G)$}{spanning subgraph of $G$ with
  edge set $E$} %
\nomenclature[4g2e]{$G-e, e\in E(G)$}{spanning subgraph of $G$ with edge
  set $E(G)\setminus \{e\}$} %

By an {\em induced subgraph}, we always mean a subgraph induced by a
vertex set; a subgraph induced by an edge set is called an {\em
  edge-subgraph}. %
We write $H \subseteq G$ when $H$ is a subgraph of $G$, %
\nomenclature[1as]{$\subseteq$}{subset of; subgraph of} %
$H \subseteq_e G$ when $H$ is an edge-subgraph of $G$, %
\nomenclature[1ase]{$\subseteq_e$}{edge-subgraph of} %
$H \subseteq_v G$ when $H$ is an induced subgraph of $G$, %
\nomenclature[1asv]{$\subseteq_v$}{induced subgraph of} %
$H \leq G$ when $H$ is isomorphic to a subgraph of $G$, %
\nomenclature[1al]{$\leq$}{less than or equal; isomorphic to a subgraph
  of} %
$H \epo G$ when $H$ is isomorphic to an edge-subgraph of $G$, and %
\nomenclature[1ale]{$\epo$}{isomorphic to an edge-subgraph of} %
$H \vpo G$ when $H$ is isomorphic to an induced subgraph of $G$. %
\nomenclature[1alv]{$\vpo$}{isomorphic to an induced subgraph of} %
We denote the number of subgraphs (induced subgraphs, edge-subgraphs,
components) of $G$ that are isomorphic to $H$ by $\sub{H}{G}$
(respectively, $\vsub{H}{G}$, $\esub{H}{G}$, $c(H,G)$). %
\nomenclature[4sub]{$\sub{H}{G}$}{number of subgraphs of $G$ that are
  isomorphic to $H$} %
\nomenclature[4i3nd]{$\vsub{H}{G}$}{number of induced subgraphs of $G$
  that are isomorphic to $H$} %
\nomenclature[4c3hg]{$c(H,G)$}{number of components of $G$ that are
  isomorphic to $H$} %

\subsubsection{Unlabelled graphs}
If $H$ is isomorphic to $G$, then we write $H\cong G$. %
\nomenclature[1]{$\cong$}{isomorphic to; used for graphs and posets} %
Isomorphism is an equivalence relation on $\lab$. An {\em unlabelled
  graph} is an isomorphism class. A {\em class of graphs} is a set of
graphs closed under isomorphism.
We use the quotient notation $\unlab$ to denote the set of {\em
  isomorphism classes}, %
\nomenclature[4g1]{$\unlab $}{set of all unlabelled graphs} %
and take $I$ to be the quotient map of isomorphism; thus
$I(G) \coloneqq \{H \in \lab\mid H \cong G\}$ is the isomorphism class
of $G$. For $S \subseteq \lab$, we write
$\iso{S} \coloneqq \{\iso{G} \mid G \in S\}$. %
Given an isomorphism class (unlabelled graph) $G$, we denote a
representative labelled graph in $G$ by $\rep{G}$.
\nomenclature[4i2]{$\iso{G}$, $G\in \lab$}{isomorphism class of a
  labelled graph $G$} %
\nomenclature[4g2]{$\rep{G}$, $G\in \unlab$}{representative labelled
  graph in an isomorphism class $G$}

A {\em graph-invariant} is a function $f$ on $\lab$ that is constant
over each isomorphism class. If $f$ is a graph invariant, we define
$f(S) \coloneqq f(\rep{S})$ for all $S \in \unlab$.  Definitions of many
terms (e.g., deck, edge-deck, etc.)  and parameters (e.g.,
$\vsub{.}{.}$, $\esub{.}{.}$, $\hom(.,.)$, etc.)  naturally extend to
and are well-defined for unlabelled graphs if they depend on invariant
properties of graphs. For example, if $S,T,\in \lab$, then $\vsub{S}{T}
= \vsub{S^\prime}{T^\prime}$ for all $S^\prime\in \iso{S}$, for all
$T^\prime\in \iso{T}$, which allows us to define $\vsub{S}{T}\coloneqq
\vsub{\rep{S}}{\rep{T}}$ if $S$ and $T$ are unlabelled
graphs. Similarly, for unlabelled graphs $S$ and $T$, we say that $S$ is
an induced subgraph (or an edge-subgraph) of $T$ if $\rep{S}$ is
isomorphic to an induced subgraph (or an edge-subgraph) of $\rep{T}$.

We denote a path on $n$ vertices by $P_n$, %
a cycle on $n$ vertices by $C_n$, %
a complete graph on $n$ vertices by $K_n$, %
a complete bipartite graph with $n$ and $m$ vertices in the two
partitions by $K_{n,m}$, %
and the graph $K_4$ minus an edge by $K_4\setminus e$; %
\nomenclature[4p2n]{$P_n$}{path on $n$ vertices} %
\nomenclature[4c2n]{$C_n$}{cycle on $n$ vertices} %
\nomenclature[4k2n]{$K_n$}{complete graph on $n$ vertices} %
\nomenclature[4k2n]{$K_{n,m}$}{complete bipartite graph} %
\nomenclature[4k2e]{$\eminus{K_4}$}{$K_4$ minus an edge} %
here $P_k$, $C_k$ and $K_4\setminus e$ are unlabelled graphs. Similarly,
$K_n$, $K_{n,m}$ are unlabelled graphs. We write $G \in K_n$
\mc{labelled} to refer to a (labelled) graph in $K_n$.

Let $H_i, i =\seq{1}{m}$ be distinct unlabelled graphs.  Let
$\mcf \coloneqq \{F_1,\dotsc,F_n\} \subseteq \lab$ be a collection of
mutually vertex-disjoint graphs. If $G \coloneqq \biguplus_{i=1}^n F_i$
and $k_i \coloneqq \card{\mcf\cap H_i}$, for $i=\seq{1}{m}$, then we
write $G\in \sum_i k_iH_i$ and $\iso{G} = \sum_i k_iH_i$.

\subsubsection{Reconstruction terminology} 
Most of the following notions are standard in the reconstruction
literature (see, e.g., \citet{bondy1991}), so we define them concisely
below.

The deck of a labelled graph $G$ is the set
$\deck{G}\coloneqq \{(\iso{G-u}, \vsub{G-u}{G}) \mid u \in V(G)\}$. We
write $\deck{G}\coloneqq \deck{\rep{G}}$ when $G$ is an unlabelled
graph. We say that $H$ is a reconstruction of $G$ if
$\deck{G} = \deck{H}$; and that $G$ is {\em reconstructible} if it is
determined up to isomorphism by $\deck{G}$ (i.e., every reconstruction
of $G$ is isomorphic to $G$).  A set $S$ of unlabelled graphs is a {\em
  counter example} to Ulam's conjecture if $\deck{G} = \deck{H}$ for all
$G,H \in S$ and $\card{S} \geq 2$. %
Let $\mcc$ be a class of graphs and let $f$ be a graph invariant. We say
that %
$\mcc$ is reconstructible if each graph in $\mcc$ is reconstructible; %
$f(G)$ is reconstructible if $\deck{G}$ determines $f(G)$; and %
$f$ is reconstructible for $\mcc$ if it is reconstructible for all
graphs in $\mcc$.

Similar definitions may be given for other reconstruction problems. In
particular, by replacing {\em deck} by {\em edge-deck} or by {\em
  abstract induced subgraph poset} (to be defined in
Section~\ref{sec-intro-isp}) or by {\em abstract bond lattice } (to be
defined in Section~\ref{sec-intro-lattice}) or by {\em abstract
  edge-subgraph poset} (to be defined in Section~\ref{sec-intro-esp}),
we define the corresponding notions of {\em edge reconstructibility},
{\em $P$-reconstructibility}, {\em $\Pi$-reconstructibility}, {\em
  $Q$-reconstructibility}, respectively.

\subsubsection{Partially ordered sets}
\label{sec:posets}
We follow \citet{stanley-v1} for terminology on partially ordered
sets. Let $(S,\leq)$ be a partially ordered set. For $x, y \in S$, we
say that $y$ {\em covers} $x$ if $x \leq y$, $x \neq y$, and there is no
$z\in S\setminus \{x,y\}$ such that $x \leq z \leq y$. A partially
ordered set (poset) is called {\em ranked} if it admits a {\em rank
  function} $\rho : S\to \nn$ such that for all $x,y$ in $S$, $y$ covers
$x$ implies $\rho(y) = \rho(x) + 1$. %
The {\em down-set} of an element $x$ is the set
$\down{S}{x} \coloneqq \{y \in S \colon\, y \leq x\}$.

The posets in this paper are weighted. We say that weighted posets
$(S, \leq, w)$ and $(S^\prime, \leq^\prime, w^\prime)$ are {\em
  isomorphic} if there is a bijection $f\colon S \to S^\prime$, called
an {\em isomorphism}, such that for all $x,y\in S$, we have $x \leq y$
if and only if $f(x) \leq^\prime f(y)$ and
$w(x,y) = w^\prime(f(x),f(y))$. An {\em automorphism} of a poset
$(S, \leq, w)$ is an isomorphism from $(S, \leq, w)$ to itself.

\subsubsection{Partitions}
\label{sec:partitions}
A partition of a positive integer $n$ is a tuple
$\lambda \coloneqq \tup{\lambda_1}{\lambda_k}$ of integers, where
$\lambda_1 \geq \cdots \geq \lambda_k > 0$, and $\sum_i\lambda_i = n$.
The length of a partition $\lambda$, denoted by $\ell(\lambda)$, is the
number of elements in $\lambda$. %
\nomenclature[4l]{$\ell({\lambda})$}{length of an integer partition} %
Let $\lambda \coloneqq \tup{\lambda_1}{\lambda_k}$ and
$\mu \coloneqq \tup{\mu_1}{\mu_l}$ be two partitions of $n$. We say that
$\lambda$ {\em refines} $\mu$, and write $\lambda \models \mu$, %
\nomenclature[5lm]{$\lambda \models \mu$}{integer partition $\lambda$
  refines integer partition $\mu$} %
if there is an onto map $f\colon \ints{1}{k}\to \ints{1}{l}$ such that
$\sum_{j \in f^{-1}(i)} \lambda_j = \mu_i$ for all $i\in \ints{1}{l}$.
The refinement relation makes the set of partitions of $n$ a lattice. We
also write $\lambda\models n$ to say that $\lambda$ is a partition of
$n$. %
\nomenclature[5ln]{$\lambda \models n$}{$\lambda$ is a partition of
  $n$} %

Let $\pi \coloneqq \set{X_1}{X_k}$ be a family of mutually disjoint
non-empty subsets of $V$. We write $\pi \vdash V$. %
\nomenclature[5pi2v2]{$\pi \vdash V$}{$\pi$ is a family of mutually
  disjoint subsets of $V$} %
If $\pi$ is a partition of $V$, we write $\pi \models V$; %
\nomenclature[5pi2v2]{$\pi \models V$}{$\pi$ is partition of $V$} %
we associate with $\pi$ an integer partition $\lambda(\pi)$ of $|V|$
obtained by ordering $|X_i|$ in a non-increasing order. %
\nomenclature[5lp]{$\lambda(\pi)$}{integer partition associated with a
  set partition $\pi$} %
Let $\pi$ and $\sigma$ be partitions of $V$.  We say that $\pi$ {\em
  refines} $\sigma$, and write $\pi \models \sigma$, if each block of
$\pi$ is a subset of some block of $\sigma$. %
\nomenclature[5pi2sigma]{$\pi \models \sigma$}{set partition $\pi$
  refines set partition $\sigma$} %
The refinement relation makes the set of partitions of $V$ a lattice. It
is called the {\em partition lattice} of $V$, and is denoted by
$\Pi(V)$. %
\nomenclature[5pi1v1]{$\Pi(V)$}{partition lattice of $V$}

\subsection{Induced subgraph posets and Ulam's conjecture}
\label{sec-intro-isp}

\begin{defn}
  Define a partial order $\vpo$ on $\unlab$ as follows: for all
  $F_i,F_j \in \unlab$, $F_i \vpo F_j$ if and only if $F_i$ is an
  induced subgraph of $F_j$.
  Define a {\em weight function}
  $\wv\colon(\unlab)\times(\unlab) \to \nn$ as follows: for all
  $F_i,F_j \in \unlab$, $\wv(F_i,F_j) \coloneqq \vsub{F_i}{F_j}$.
  Let $\ispug\coloneqq (\unlab,\vpo,\wv)$.
  \nomenclature[4p1g]{$\ispug$}{induced subgraph poset with ground set
    $\unlab$}
  \nomenclature[4wv]{$\wv$}{weight function of the inducted subgraph poset} %
  For $G \in \unlab$, let
  $\ispu{G}\coloneqq\{F\in\unlab\mid F\vpo G, \text{ and } \nu(F)=1
  \text{ or } \epsilon(F)>0\}$. %
  The {\em concrete induced subgraph poset} of $G$ is the restriction of
  $\ispug $ to $\ispu{G}$; it is denoted by just $\ispu{G}$. %
  \nomenclature[4p1g]{$\ispu{G}$}{concrete induced subgraph poset of
    $G$} The {\em abstract induced subgraph poset} is the isomorphism
  class of $(\ispu{G}, \vpo, \wv)$. %
  We take $\ispu{G} \coloneqq \set{G_1}{G_M}$, where $\seq{G_1}{G_M}$
  are distinct unlabelled graphs, and
  $\isp{G} \coloneqq (\set{g_1}{g_{\scriptscriptstyle M}},\vpo,\wv)$ to
  be a representative poset isomorphic to $\ispu{G}$. %
  \nomenclature[4p1g]{$\isp{G}$}{abstract induced subgraph poset of
    $G$} %
  \nomenclature[4g3]{$g_1,\ldots,g_{\scriptscriptstyle M}$}{elements in
    an abstract induced subgraph poset} %
  Moreover, we assume that there is an isomorphism
  $f\colon \isp{G} \to \ispu{G}$ such that $f(g_i) = G_i$ for all $i$.
  We assume that the minimum element in $\isp{G}$ is $g_1$ (thus
  $G_1 = K_1$), and that $g_2$ covers $g_1$ (thus $G_2 = K_2$). Since
  $\epsilon(G_i) = \wv(g_2,g_i)$, we define
  $\epsilon(g_i) \coloneqq \wv(g_2,g_i)$. We assume that the maximal
  element in $\isp{G}$ is $g_{\scriptscriptstyle M}$ (thus $G_{M} =
  G$).
  We define a rank function $\nu$ on $\isp{G}$ so that
  $\nu(g_1) \coloneqq 1$; thus $\nu(g_i)=\nu(G_i)$ for all $i$. %
\end{defn}

\begin{rem} \mc{labelled} When a graph $G$ is labelled, we define
  $\ispu{G} \coloneqq \ispu{\iso{G}}$, and
  $\isp{G} \coloneqq \isp{\iso{G}}$, and so on, where $\iso{G}$ is the
  isomorphism class of $G$. The set $\ispu{G}$ contains only unlabelled
  graphs.
\end{rem}

\begin{defn}
  An unlabelled graph $G$ is said to be {\em $P$-reconstructible} if it
  is determined by $\isp{G}$. A labelled graph $G$ is said to be {\em
    $P$-reconstructible} if its isomorphism class is determined by
  $\isp{G}$. A class of graphs is said to be {\em $P$-reconstructible}
  if each graph in the class is $P$-reconstructible.
\end{defn}

\begin{figure}[ht]
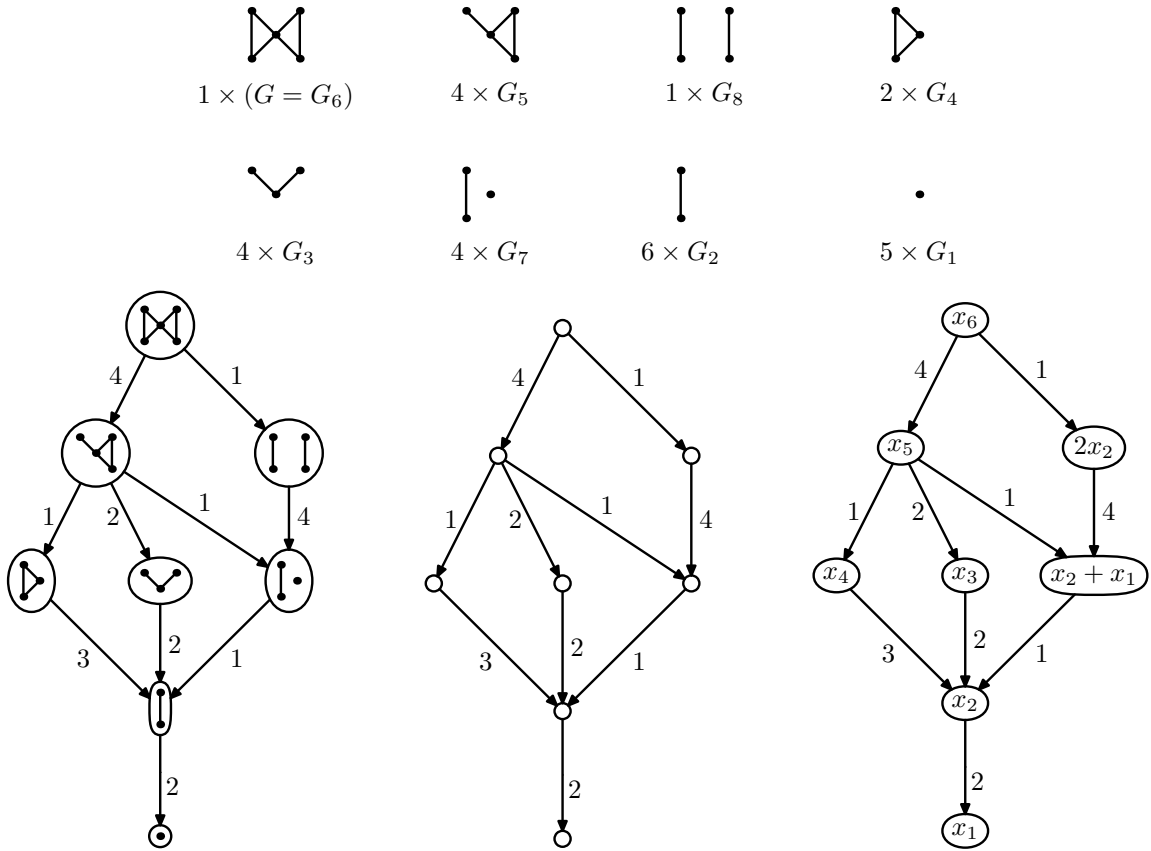

\begin{center}
\includegraphics{figs-1.mps}

\vspace{0.1in}
\includegraphics{figs-2.mps}
\hspace{0.4in}
\includegraphics{figs-3.mps}
\hspace{0.4in}
\includegraphics{figs-4.mps}

\end{center}
\caption[]{An unlabelled graph $G$ and the graphs in its concrete
  induced subgraph poset along with their multiplicities are shown in
  the top two rows. In the bottom row, the concrete induced subgraph
  poset is shown on the left; the abstract induced subgraph poset is
  shown in the middle; a consistent labelling of its elements by formal
  sums of indeterminates (see Definition~\ref{defn-labels}) is shown on
  the right.}
\label{fig-pg}
\end{figure}

\begin{example}
  An unlabelled graph $G$, the graphs in its concrete induced subgraph
  poset along with their multiplicities, and the induced subgraph poset
  of $G$ (concrete and abstract) are shown in Figure~\ref{fig-pg}. Note
  that the Hasse diagram is only for illustration; it does not display
  all weights. But the weights on all related pairs can be computed from
  the weights shown in the Hasse diagram; see Lemma 2.3 in
  \cite{thatte2005}.
\end{example}

In \cite{thatte2005}, we considered the {\em $P$-reconstruction
  problem}, i.e., the problem of reconstructing $G$ up to isomorphism or
computing some of its invariants from $\isp{G}$. The following theorem
summarises our results from \cite{thatte2005}.

\begin{thm}
  \label{thm-thatte} \text{}
  \begin{enumerate}
  \item\label{thm-thatte-ulam} Ulam's conjecture is true if and only if
    every non-empty graph can be reconstructed up to isomorphism from
    its abstract induced subgraph poset.
  \item\label{thm-thatte-auto} Ulam's conjecture is true if and only if
    the abstract induced subgraph poset of every non-empty graph has only the
    trivial automorphism.
  \item\label{thm-thatte-trees} Every tree can be reconstructed up to
    isomorphism from its abstract induced subgraph poset.
  \item\label{thm-thatte-invariants} The following invariants of a graph
    $G$ can be computed from its abstract induced subgraph poset:
    \begin{enumerate}
    \item the number of spanning trees in $G$; and hence whether
      $G$ is connected,
    \item the number of unicyclic subgraphs of $G$
      containing a cycle of a specified length; and hence the number of
      Hamiltonian cycles in $G$,
    \item the number of spanning subgraphs of $G$ having specified
      numbers of vertices and edges in their components,
    \item the characteristic polynomial of $G$, the chromatic
      polynomial of $G$, and the rank polynomial of $G$.
    \end{enumerate}
  \end{enumerate}
\end{thm}

The invariants listed above were originally proved to be reconstructible
by \citet{tutte1979}. The above results are slightly stronger than the
results of Tutte in the sense that to compute the invariants listed
above for a graph, we do not need to know its deck - its abstract
induced subgraph poset is sufficient.

\subsection{The connected partition lattice}
\label{sec-intro-lattice}

Computing the graph invariants listed in
Theorem~\ref{thm-thatte}-(\ref{thm-thatte-invariants}) (in our proofs as
well as in the proofs by \citet{tutte1979} and \citet{kocay1981})
requires first counting certain disconnected spanning subgraphs. When a
deck is given, counting disconnected spanning subgraphs is made easier
by a lemma of \citet{kocay1981}. In the proof of
Theorem~\ref{thm-thatte}, since the deck is not given, counting
disconnected spanning subgraphs with a given number of components, and a
given number of vertices and edges in each component, is more
difficult. We can nevertheless imitate Kocay's lemma by using minimal
information about the graphs in the induced subgraph poset (e.g.,
$\nu(G_i)$ and $\epsilon(G_i)$ for each $G_i$ in $\ispu{G}$). Counting
disconnected spanning subgraphs would be easier if we were given the
lattice of connected partitions.

In the proof of Theorem~\ref{thm-thatte}-(\ref{thm-thatte-invariants}),
we implicitly constructed and used partial information about the
connected partition lattice. We commented in \cite[Section
5]{thatte2005} that understanding the relationship between the induced
subgraph poset of a graph and its connected partition lattice would be
interesting. One of the objectives of this paper is to clarify this
relationship.

First we define the connected partition lattice of a labelled
graph. Recall that when $\pi\coloneqq\set{X_1}{X_m}$ is a family of
disjoint non-empty subsets of a set $V$, we write $\pi\vdash V$; and
when $\pi\coloneqq\set{X_1}{X_m}$ is a partition of $V$, we write
$\pi \models V$.
\begin{defn}
  \mc{labelled} Let $G$ be a labelled graph. Let
  $\pi\coloneqq\set{X_1}{X_m}$ be a family of disjoint non-empty subsets
  of $V(G)$. We say that $\pi$ is a {\em connected family} (or {\em
    $\pi$ is a connected partition} in case $\pi$ is a partition of
  $V(G)$) if subgraphs $G[X_k], k \in [m]$ induced by the blocks of
  $\pi$ are all connected; in this case we write $\pi \latpo V(G)$ (or
  $\pi \models_c V(G)$ in case $\pi$ is a partition of $V(G)$). We
  define $G[\pi]\coloneqq\biguplus_{k\in[m]}G[X_k]$, and call it the
  {\em subgraph of $G$ induced by $\pi$}.
  Let $\latpo$ be a partial order on $\lab$ defined by: $H \latpo G$ if
  there exists a connected family $\pi \latpo V(G)$ such that
  $G[\pi] = H$.
  Let $\latl{G}$ denote the set of subgraphs of $G$ induced by its
  connected partitions, i.e.,
  $\latl{G}\coloneqq \{G[\pi]\mid \pi \modpo V(G)\}$. The {\em connected
    partition lattice} of $G$ is the restriction of $(\lab,\latpo)$ to
  $\latl{G}$; it is denoted by just $\latl{G}$. %
  \nomenclature[5pi1g]{$\latl{G}$}{connected partition lattice of $G$} %
  The unique minimal element $\hat{0}$ of $\latl{G}$ is the finest
  partition of $V(G)$ consisting of only singletons. %
  \nomenclature[20]{$\hat{0}$}{minimal element of the connected
    partition lattice} %
\end{defn}

Next we define analogous notions for an unlabelled graph. Recall that
given an unlabelled graph $H$, we denote a representative labelled graph
in $H$ by $\rep{H}$.

\begin{defn}
  The partial order $\latpo$ on the set of labelled graphs induces a
  partial order $\latpo$ on the set of unlabelled graphs. We define it
  by considering labelled representatives of unlabelled graphs.
  For unlabelled graphs $H_j, H_k$, we define $H_j \latpo H_k$ if there
  exists $\pi \latpo V(\rep{H_k})$ such that
  $\rep{H_k}[\pi] \cong \rep{H_j}$. %
  We define a {\em weight function}
  $\wl\colon (\unlab)\times (\unlab) \to \nn$ by
  $\wl(H_j,H_k) \coloneqq |\{\pi\latpo V(\rep{H_k}) \mid \rep{H_k}[\pi]
  \in H_j\}|$.
  Thus we have a weighted partially ordered set
  $\latug \coloneqq (\unlab,\latpo, \wl)$.
  The {\em folded connected partition lattice} of a labelled graph $G$
  is the restriction of $\latug $ to $\iso{\latl{G}}$. %
  (Recall that $\iso{\latl{G}}$ is the set of distinct isomorphism
  classes of graphs in $\latl{G}$.) %
  We denote the set $\iso{\latl{G}}$ as well as the folded connected
  partition lattice of $G$ by just $\latu{G}$. %
  \nomenclature[5z3g]{$\latu{G}$}{concrete bond lattice of $G$} %
  \nomenclature[4wl]{$\wl$}{weight function of the bond lattice} %
  The {\em abstract folded connected partition lattice} of $G$ is the
  isomorphism class of $\latu{G}$. %
  \nomenclature[5z3g]{$\lat{G}$}{abstract bond lattice of $G$} %
  We take $\latu{G} \coloneqq \set{H_1}{H_N}$, where $\seq{H_1}{H_N}$ are
  distinct unlabelled graphs, and
  $\lat{G} \coloneqq (\set{h_1}{h_{\scriptscriptstyle N}},\modpo,\wl)$
  to be a representative poset isomorphic to $\latu{G}$. Moreover, we
  assume that there is an isomorphism $f\colon \lat{G} \to \latu{G}$
  such that $f(h_i) = H_i$ for all $i$.
  We assume that the minimal element in $\lat{G}$ is $h_1$ (thus
  $H_1 = \nu(G)K_1$), and the maximal element in $\lat{G}$ is
  $h_{\scriptscriptstyle N}$ (thus $H_{N} = \iso{G}$).  We define a rank
  function $\rho$ on $\lat{G}$ so that $\rho(h_1) = 0$; hence
  $c(H_i) = \nu(G) - \rho(h_i)$ for all $i$. %
  \nomenclature[5rho]{$\rho(x)$}{rank of element $x$ in a poset} %
  \nomenclature[4h]{$h_1,\ldots,h_N$}{elements in an abstract bond
    lattice} %
\end{defn}

\begin{rem} \mc{labelled} When $G$ is an unlabelled graph, we define
  $\latu{G}$ and $\lat{G}$ in terms of a representative labelled graph;
  i.e., $\latu{G} \coloneqq \latu{\rep{G}}$ and
  $\lat{G} \coloneqq \lat{\rep{G}}$. Regardless of whether $G$ is
  labelled or unlabelled, the set $\latu{G}$ contains only unlabelled
  graphs.
\end{rem}

\begin{defn}
  An unlabelled graph $G$ is said to be {\em $\Pi$-reconstructible} if
  it is determined by $\lat{G}$. A labelled graph $G$ is said to be {\em
    $\Pi$-reconstructible} if its isomorphism class is determined by
  $\lat{G}$. A class of graphs is said to be {\em $\Pi$-reconstructible}
  if each graph in the class is $\Pi$-reconstructible.
\end{defn}

\begin{rem}
  The connected partition lattice of a graph $G$ forms a geometric
  sub-lattice of the partition lattice $\Pi(V(G))$. Elsewhere in the
  literature, it has been called the {\em lattice of contractions} or
  the {\em bond lattice} of $G$; see
  \citet{stanley1995,stanley-v1,chow1995}. The lattice $\latu{G}$
  (without its weights) is obtained by identifying those partitions in
  $\latl{G}$ that induce isomorphic graphs; hence we think of $\latu{G}$
  as the {\em folded connected partition lattice} of $G$. We use the
  term {\em abstract bond lattice} only for the abstract folded
  connected partition lattice, and the term {\em concrete connected
    partition lattice} for the labelled lattice $\latl{G}$ defined
  above.
\end{rem}

\begin{figure}[ht]
\begin{center}
\includegraphics{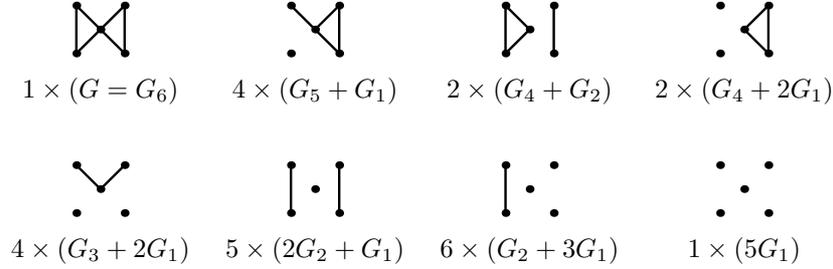}
\end{center}
\caption[]{An unlabelled graph $G$ and its distinct connected partitions
  along with their multiplicities.}
\label{fig-pisub}
\end{figure}

\begin{figure}[ht]
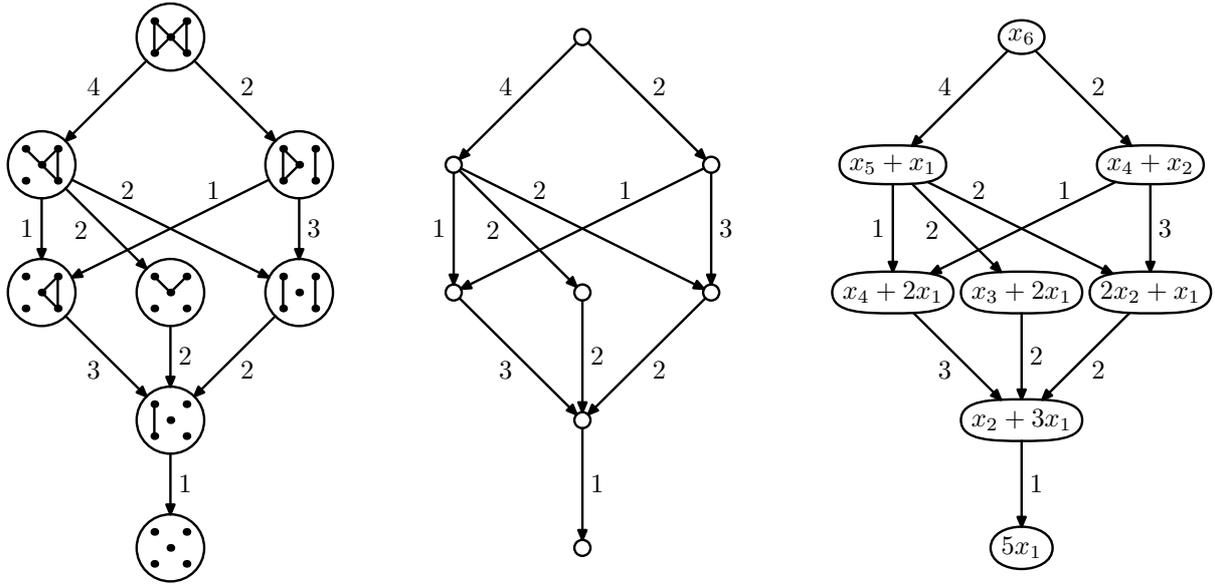

\begin{center}
\includegraphics{figs-6.mps}
\hspace{0.4in}
\includegraphics{figs-7.mps}
\hspace{0.4in}
\includegraphics{figs-8.mps}
\end{center}
\caption[]{The connected partition lattice of an unlabelled graph $G$:
  the lattice on the left is concrete; the lattice in the middle is
  abstract; a consistent labelling of its elements with formal sums of
  indeterminates (see Definition~\ref{defn-labels}) is shown on the
  right.}
\label{fig-lc}
\end{figure}

\begin{example}
  An unlabelled graph $G$ and its connected partitions along with their
  multiplicities are shown in Figure~\ref{fig-pisub}. The graphs
  $G_1,G_2,\dotsc,G_6=G$ are the same graphs as in
  Figure~\ref{fig-pg}. %
  The bond lattice of $G$ (at different levels of abstraction) is shown
  in Figure ~\ref{fig-lc}. %
  To be precise, the lattice in the middle is abstract folded connected
  partition lattice or the abstract bond lattice. %
  A consistent labelling of its elements by formal sums of
  indeterminates, indicating the component structure of graphs induced
  by connected partitions, is shown on the right. The notion of
  consistent labelling is formalised in Definition~\ref{defn-labels}. %
  As in the case of the induced subgraph poset, the Hasse diagram is
  only for illustration; it does not display the weights on all related
  pairs of graphs in $\lat{G}$.
\end{example}

\begin{example} \label{ex-nk2-k1n} The graphs $K_{1,n}$ and $nK_2$ have
  the same abstract bond lattice, which is a total order with
  appropriate weights. Moreover, adding isolated vertices to a graph
  does not change its abstract bond lattice. But no non-empty graph
  other than graphs isomorphic to $nK_2 + mK_1$ or $K_{1,n} + mK_1$ for
  some $n,m \in \nn$ has a totally ordered abstract bond lattice. We do
  not know any other non-trivial pairs of nonisomorphic graphs that have
  the same abstract bond lattice. On the other hand, the graphs
  $K_{1,n}$ and $nK_2$ have distinct induced subgraph posets.
\end{example}

In Section~\ref{sec-kocay}, we introduce Kocay's lemma
(Lemma~\ref{lem-kocay1}) and two variants of it, namely,
Lemma~\ref{lem-kocay2} for induced subgraphs and Lemma~\ref{lem-kocay3}
for abstract induced subgraph posets. Given a list of graphs
\mc{labelled} $\tup{F_1}{F_k} \in \lab^k$ and the deck of $G$, Kocay's
lemma (see \citet{bondy1991} or \citet{kocay1981}) shows how to compute
the number of {\em covers of $G$ by $\tup{F_1}{F_k}$}, i.e., the number
of tuples $\tup{H_1}{H_k}$ of subgraphs of $G$ such that for all
$i=\seq{1}{k}$ we have $H_i \cong F_i$ and $\bigcup_i H_i = G$.

Our generalisation in Lemma~\ref{lem-kocay3} demonstrates how we may
compute the number of covers by a list of unknown induced subgraphs of
$G$ (i.e., the graphs $\iso{F_i}$ are specified only as elements of the
abstract induced subgraph poset). We illustrate Lemma~\ref{lem-kocay3}
by giving simple proofs that for every nonempty graph $G$, the chromatic
symmetric function of $G$ and the symmetric Tutte polynomial of $G$
(both introduced by \citet{stanley1995,stanley1998}) are
$P$-reconstructible.

In Section~\ref{sec-pg-lg}, we give two constructions: we show, using
counting arguments based on Kocay's lemma that are developed in
Section~\ref{sec-kocay}, that the abstract bond lattice of every graph
can be constructed from its abstract induced subgraph poset
(Theorem~\ref{thm-pg2lg}); and that the abstract induced subgraph poset
of every graph that is not a star or a disjoint of edges and that has no
isolated vertices can be constructed from its abstract bond lattice
(Theorem~\ref{thm-lg2pg}). Theorem~\ref{thm-pg2lg} and
Theorem~\ref{thm-thatte}-(\ref{thm-thatte-trees}) together imply
Corollary~\ref{cor-ulam}-(\ref{cor-ulam-trees}) that every tree (or
forest) on 2 or more vertices that is not a star or a disjoint union of
edges and that has no isolated vertices can be reconstructed up to
isomorphism from its abstract bond lattice. In
Section~\ref{sec-polynomials}, we give short proofs of the
reconstructibility of the symmetric Tutte polynomial and the chromatic
symmetric function. In particular, we give another expansion of the
chromatic symmetric function based on the abstract bond lattice
(Corollary~\ref{cor-xg3}).

Another motivation of this paper is a question of Stanley regarding the
chromatic symmetric function: can two non-isomorphic trees have the same
chromatic symmetric function?
Corollary~\ref{cor-ulam}-(\ref{cor-ulam-trees}) suggests that one way to
approach Stanley's question may be to try to construct the abstract bond
lattice of a tree from its chromatic symmetric function.  In
Section~\ref{sec-trees}, we give a few preliminary results, which we
summarise below.

Let $T$ be a tree with vertex set $V(T)$. Let
$\mathbf{v} \coloneqq \tup{v_1}{v_r}$ and
$\mathbf{e} \coloneqq \tup{e_1}{e_r}$ be two integer vectors. Let
$\theta(\mathbf{v},\mathbf{e};T)$ be the number of ordered partitions
$\tup{V_1}{V_r}$ of $V(T)$ such that $|V_i| = v_i$, and
$\epsilon(T[V_i]) = e_i$. We show in Lemma~\ref{lem-theta} that for all
trees $T$ and for all vectors $\mathbf{v}$ and $\mathbf{e}$, we can
compute $\theta(\mathbf{v},\mathbf{e};T)$ from the chromatic symmetric
function $X_T(x)$ of $T$. We apply the lemma to prove that the symmetric
Tutte polynomial of every tree can be obtained from its chromatic
symmetric function, which is not the case for graphs in
general. \citet{noble_welsh1999} have asked if trees are distinguished
by their U-polynomial, which, for trees, is equivalent to the symmetric
Tutte polynomial. Our result shows that the questions of Noble \& Welsh
and Stanley are in fact equivalent. The collection of invariants
$\theta(\mathbf{v},\mathbf{e};T)$ may have sufficient information to
construct the abstract bond lattice, and hence the tree by
Corollary~\ref{cor-ulam}-(\ref{cor-ulam-trees}).

\subsection{The edge-subgraph poset}
\label{sec-intro-esp}

\begin{defn}
  Define a partial order $\epo$ on $\unlab$ as follows: for all
  $F_i,F_j \in \unlab$, $F_i \epo F_j$ if and only if $F_i$ is an
  edge-subgraph of $F_j$.
  Define a {\em weight function} $\we \colon (\unlab)\times(\unlab) \to \nn$
  as follows: for all $F_i,F_j \in \unlab$,
  $\we(F_i,F_j) \coloneqq \esub{F_i}{F_j}$.
  Let $\espug\coloneqq (\unlab,\epo,\we)$. %
  \nomenclature[4qg]{$\espug$}{weighted edge-subgraph poset with ground
    set $\unlab$} %
  \nomenclature[4we]{$\we$}{weight function of the edge-subgraph
    poset} %

  For $G \in \unlab$, let
  $\espu{G} \coloneqq \{F \in \unlab \mid F \epo G \text{ and }
  \epsilon(F) > 0\}$.
  The {\em concrete edge-subgraph poset} of $G$ is the restriction of
  $\espug$ to $\espu{G}$; it is denoted by just $\espu{G}$. %
  \nomenclature[4qg]{$\espu{G}$}{concrete edge-subgraph poset} %
  The {\em abstract edge-subgraph poset} of $G$ is the isomorphism class
  of $(\espu{G},\epo,\we)$. %
  We take $\espu{G} \coloneqq \set{G_1}{G_M}$, where $\seq{G_1}{G_M}$
  are distinct unlabelled graphs, and
  $\esp{G} \coloneqq (\set{g_1}{g_{\scriptscriptstyle M}}, \epo, \we)$
  to be a representative poset isomorphic to $\espu{G}$. %
  \nomenclature[4qg]{$\esp{G}$}{abstract edge-subgraph poset} %
  We assume that there is an isomorphism $f\colon \esp{G} \to \espu{G}$
  such that $f(g_i) = G_i$ for all $i$; hence
  $\esp{g_i} \cong \espu{G_i}$ for $i = 1\text{ to } M$.
  We assume that the minimum element in $\esp{G}$ is $g_1$ (thus
  $G_1 = K_2$), and that the maximal element in $\esp{G}$ is
  $g_{\scriptscriptstyle M}$ (thus $G_{M} = G$). We define a rank
  function $\rho$ on $\esp{G}$ so that $\rho(g_1) \coloneqq 1$; thus
  $\rho(g_i)=\epsilon(G_i)$ for all $i$.
  \nomenclature[4g3]{$g_1,\ldots,g_{\scriptscriptstyle M}$}{elements in
    an abstract edge-subgraph poset} %
\end{defn}

\begin{rem} \mc{labelled} When a graph $G$ is labelled, we define
  $\espu{G} \coloneqq \espu{\iso{G}}$, and
  $\esp{G} \coloneqq \esp{\iso{G}}$, where $\iso{G}$ is the isomorphism
  class of $G$. The set $\espu{G}$ contains only unlabelled graphs.
\end{rem}

\begin{defn}
  An unlabelled graph $G$ is said to be {\em $Q$-reconstructible} if it
  is determined by $\esp{G}$. A labelled graph $G$ is said to be {\em
    $Q$-reconstructible} if its isomorphism class is determined by
  $\esp{G}$. A class of graphs is said to be {\em $Q$-reconstructible}
  if each graph in the class is $Q$-reconstructible. %
  For all $G \in \unlab$, we have $\esp{G} = \esp{G+K_1}$.  Therefore,
  we understand $Q$-reconstructibility to mean $Q$-reconstructibility
  {\em modulo isolated vertices}. %
  A {\em $Q$-set} is a maximal set of cardinality at least 2 of
  unlabelled graphs that have no isolated vertices, and have isomorphic
  edge-subgraph poset. A {\em $Q$-pair} is 2-element subset of a
  $Q$-set. %
\end{defn}

We ask which graphs are $Q$-reconstructible, and if there is a relation
between the problem of $Q$-reconstructibility and the edge
reconstruction conjecture of Harary (analogous to
Theorem~\ref{thm-thatte}-(\ref{thm-thatte-invariants}) for Ulam's
conjecture). %
The question of $Q$-reconstructibility is not quite similar to the
analogous question for induced subgraphs. If a graph is not edge
reconstructible, then it is also not $Q$-reconstructible. It turns out
that there many graphs that are edge reconstructible but not
$Q$-reconstructible.

\begin{example}
  Figure~\ref{fig-part-labelled} shows part of an edge-subgraph poset:
  its elements are labelled by the corresponding graphs on the left; a
  representative isomorphic poset is shown on the right. %
  Graphs $K_{1,2}$, $2K_2$, $K_{1,3}$, and $K_3$ are not
  $Q$-reconstructible, since they are not edge reconstructible. %
  There are many more graphs that are not $Q$-reconstructible. For
  example, the down-sets of $g_4$, $g_7$ and $g_8$ are isomorphic, and
  correspond to the $Q$-set $\{3K_2, K_{1,3}, K_3\}$. %
  In general, for all $m > 1$, $\{K_{1,m}, mK_2\}$ is a $Q$-set. The
  down-sets of $mK_2$ and $K_{1,m}$ are totally ordered, with suitable
  edge weights. %
  The down-sets of $g_{10}$ and $g_{11}$ are isomorphic, and correspond
  to the $Q$-set $\{P_4+K_2, T_4\}$, where the graph $T_4$ is defined in
  Figure~\ref{fig-names}.
\end{example}

\begin{figure}[ht]
  \begin{center}
    \includegraphics[scale=0.95]{figs-10.mps} %
    \hspace{0.4in}
    \includegraphics[scale=0.95]{figs-11.mps} %
  \end{center}
  \vspace*{-2ex}
  \caption[]{Part of an edge-subgraph poset: the poset on the left is
    labelled by the graphs corresponding to its elements; the poset on
    the right is an isomorphic representative poset.}
  \label{fig-part-labelled}
\end{figure}

In Theorem~\ref{thm-esp}, which is the main result of
Section~\ref{sec-esp}, we construct an infinite family of graphs that
are not $Q$-reconstructible, and show that the edge reconstruction
conjecture is true if and only if these are the only graphs that are not
$Q$-reconstructible.

Section~\ref{sec-hom} is motivated by edge reconstruction as well as two
results of \citet{lovasz1971} known as {\em homomorphism cancellation
  laws}. Let $\hom(G,H)$ denote the number of homomorphisms from $G$ to
$H$. Lov\'asz proved that if $G$ and $H$ are any two finite simple
graphs such that $\hom(G,F)$ = $\hom(H,F)$ for all simple graphs $F$,
then $G$ and $H$ are isomorphic. Lov\'asz also proved an analogous
complimentary result (replacing $\hom(G,F)$ by $\hom(F,G)$, and
$\hom(H,F)$ by $\hom(F,H)$ in the statement given above). Analogous to
$\isp{G}$, $\lat{G}$ and $\esp{G}$, %
we may consider the isomorphism class of a weighted complete binary
relation $\mcr$ on the set of unlabelled graphs, defined by
$\mcr\coloneqq ((\unlab) \times (\unlab), \hom)$, %
which assigns each pair $(G,H)$ of unlabelled graphs a weight
$\hom(G,H)$. %
We now ask a reconstruction-type question: does the isomorphism class of
$\mcr$ determine uniquely the unlabelled graphs corresponding to its
elements?  Equivalently, does $\mcr$ have only the trivial automorphism?
We conjecture that it does not have non-trivial automorphisms; or
equivalently, if $f\colon (\unlab)\to (\unlab)$ is a bijection such that
$\hom(G,H) = \hom(f(G),f(H))$ for all $G$ and $H$ in $\unlab$, then $f$
is the identity map. Similarly, for labelled graphs, \mc{labelled} we
conjecture that if $f:\lab \to \lab$ is a bijection such that
$\hom(G,H) = \hom(f(G),f(H))$ for all $G,H \in \lab$, then $G \cong H$.
We show in Proposition~\ref{prop-hom2} that this conjecture is weaker
than the edge reconstruction conjecture. Note also that the statement of
the conjecture is analogous to
Theorem~\ref{thm-thatte}-(\ref{thm-thatte-auto}).

\section{Kocay's lemma and its generalisations}
\label{sec-kocay}
Computations of many interesting invariants of a graph require knowing
certain spanning subgraphs. For example, to compute the characteristic
polynomial of a graph, we need to know the number of subgraphs of each
isomorphism type in the class of graphs whose connected components are
cycles or edges; see \citet{biggs1993}. \citet{kocay1981} gave a
counting argument for counting some types of spanning subgraphs of a
graph, given its deck. In this section we prove variants Kocay's lemma
and develop other counting arguments that use similar ideas as in
Kocay's lemma. The main objective here is to do similar computations on
the abstract induced subgraph poset.

\begin{defn}
\label{defn-cov}
\mc{labelled} Let $H \in \mcg$ and
$\mcf \coloneqq \tup{F_1}{F_k} \in \lab^k$.  A {\em cover} of $H$ by
$\mcf$ is a tuple $\tup{H_1}{H_k}$ of subgraphs of $H$ such that
$H_i\cong F_i$ for all $i$, and $\bigcup_i H_i = H$. Let $\cov{\mcf}{H}$
denote the number of covers of $H$ by $\mcf$. %
\nomenclature[4c3ovfh]{$\cov{\mcf}{H}$}{number of covers of graph $H$ by a
  tuple $\mcf$ of graphs} %
When $\mcf$ is a tuple of unlabelled graphs or $H$ is an unlabelled
graph, we define $\cov{\mcf}{H} \coloneqq \cov{\rep{\mcf}}{\rep{H}}$,
where $\rep{\mcf} \coloneqq (\rep{F_1},\dotsc,\rep{F_k})$.
\end{defn}

\begin{lem}{\bf (Kocay's Lemma \cite{kocay1981})}
  \label{lem-kocay1}
  \mc{labelled} Let $\mcf \coloneqq \tup{F_1}{F_k} \in \lab^k$.  We have
  \begin{equation}\label{eq-kocay1a}
    \prod_{i = 1}^{k} \sub{F_i}{G} = \sum_{H \subseteq G} \cov{\mcf}{H}
    = \sum_{H \in \unlab} \cov{\mcf}{H}\sub{H}{G}.
  \end{equation}
  Equivalently, 
  \begin{equation}\label{eq-kocay1b}
    \sum_{\substack{H \in \unlab\\
        \nu(H) \,=\, \nu(G)}} \cov{\mcf}{H}\sub{H}{G}
    = \prod_{i = 1}^{k} \sub{F_i}{G} -
    \sum_{\substack{H \subseteq G \\
        \nu(H) < \nu(G)}} \cov{\mcf}{H}\sub{H}{G}.
  \end{equation}
  Moreover, if $\nu(F_i)\, <\, \nu(G)$ for all $i$, then the left side
  of Equation~(\ref{eq-kocay1b}) is reconstructible.
\end{lem}

Next we prove Kocay's lemma for induced subgraphs.
\begin{defn}
\label{defn-vcov}
\mc{labelled} Let $H \in \mcg$ and
$\mcf \coloneqq \tup{F_1}{F_k} \in \lab^k$.  A {\em vertex cover} of $H$
by $\mcf$ is a tuple $\tup{H_1}{H_k}$ of induced subgraphs of $H$ such
that $H_i\cong F_i$ for all $i$, and $\bigcup_i V(H_i) = V(H)$.  Let
$\vcov{\mcf}{H}$ denote the number of vertex covers of $H$ by $\mcf$.
\nomenclature[4c3ovfh]{$\vcov{\mcf}{H}$}{number of vertex covers of graph
  $H$ by a tuple $\mcf$ of graphs} %
When $\mcf$ is a tuple of unlabelled graphs or $H$ is an unlabelled
graph, we define $\vcov{\mcf}{H} \coloneqq \vcov{\rep{\mcf}}{\rep{H}}$,
where $\rep{\mcf} \coloneqq (\rep{F_1},\dotsc,\rep{F_k})$.
\end{defn}
Note that a vertex cover of $H$ does not necessarily cover all edges of
$H$.

\begin{lem}{\bf (Kocay's Lemma for induced subgraphs)} 
  \label{lem-kocay2}
  \mc{labelled} Let $\mcf \coloneqq \tup{F_1}{F_k} \in \lab^k$.  We have
  \begin{equation}
    \label{eq-kocay2a}
    \prod_{i = 1}^{k} \vsub{F_i}{G} 
    = \sum_{H \subseteq_v G} \vcov{\mcf}{H}
    = \sum_{H \in \unlab}\vcov{\mcf}{H}\vsub{H}{G}.
  \end{equation}
  Moreover, if $\nu(F_i)\, <\, \nu(G)$ for all $i$, then
  $\vcov{\mcf}{G}$ is reconstructible.
\end{lem}

\begin{proof}
  Equation~(\ref{eq-kocay2a}) is self-explanatory. To prove the second
  part, we write Equation~(\ref{eq-kocay2a}) in the following form:
  \begin{equation}
    \label{eq-kocay2b}
    \vcov{\mcf}{G} = \prod_{i = 1}^{k} \vsub{F_i}{G} - 
    \sum_{\substack{H \in \unlab\\
        \nu(H) < \nu(G)}}\vcov{\mcf}{H}\vsub{H}{G},
  \end{equation}
  which is similar to Equation~(\ref{eq-kocay1b}), except that now there
  is only one term on the left side.
\end{proof}

In Lemma~\ref{lem-cov-ind} and its proof, we use the following
notation. Let $S \coloneqq \biguplus_{i=1}^m S_i$ be a graph with
connected components $S_i, i=\seq{1}{m}$. Let
$\mcs \coloneqq \tup{S_1}{S_m}$. For
$\sigma\colon\ints{1}{m}\to \ints{1}{n}$ and
$k \in \sigma(\ints{1}{m})$, we write
$S_{\sigma^{-1}(k)} \coloneqq \biguplus_{i\in \sigma^{-1}(k)}S_i$ and
$\mcs_{\sigma^{-1}(k)} \coloneqq (S_i, i \in\sigma^{-1}(k))$ (the
restricted tuple). We say that two maps
$\sigma\colon\ints{1}{m}\to \ints{1}{n}$ and
$\pi\colon\ints{1}{m}\to \ints{1}{n}$ are {\em equivalent} if
$\sigma(\ints{1}{m}) = \pi(\ints{1}{m})$ and
$S_{\sigma^{-1}(i)} \cong S_{\pi^{-1}(i)}$ for all
$i \in \sigma(\ints{1}{m})$, and {\em in-equivalent} otherwise. Let
$\Sigma $ be a set of mutually in-equivalent maps such that exactly one
representative of each equivalence class of maps is in $\Sigma$.

\begin{lem}
  \label{lem-cov-ind}
  Let $S \coloneqq \biguplus_{k=1}^m S_k$ and
  $T \coloneqq \biguplus_{k=1}^n T_k$ be two graphs with connected
  components $S_i, i \in \ints{1}{m}$ and $T_i, i \in \ints{1}{n}$,
  respectively. Let $\mcs \coloneqq \tup{S_1}{S_m}$. Then
  \begin{equation}
    \label{eq-cov1}
    \vcov{\mcs}{T}
    = \sum_{\sigma\colon\ints{1}{m}\sur \ints{1}{n}}\quad\prod_{k=1}^n 
    \vcov{\mcs_{\sigma^{-1}(k)}}{T_k},
  \end{equation}
  and
  \begin{equation}
    \label{eq-ind1}
    \vsub{S}{T} = 
    \sum_{\sigma \in \Sigma}\quad \prod_{k\,\in\, \sigma(\ints{1}{m})} 
    \vsub{S_{\sigma^{-1}(k)}}{T_k}.
  \end{equation}
\end{lem}

\begin{proof} %
  Each vertex cover $\tup{H_1}{H_m}$ of $T$ by $\mcs$ corresponds to a
  unique onto map $\sigma \colon \ints{1}{m}\to \ints{1}{n}$ such that
  $\sigma(i) = j$ if and only if $H_i \subseteq_v T_j$. Thus the set of
  vertex covers is partitioned into subsets indexed by onto maps
  $\sigma \colon \ints{1}{m}\to \ints{1}{n}$. The number of vertex
  covers that correspond to a fixed onto map $\sigma$ is the product
  term in the summation in Equation~(\ref{eq-cov1}).

  Each induced subgraph $H$ isomorphic to $S$ corresponds to a unique
  $\sigma \in \Sigma$ such that for all $k \in \sigma(\ints{1}{m})$ we
  have $H[V(T_k)] = S_{\sigma^{-1}(k)}$.  Thus the set of induced
  subgraphs of $T$ that are isomorphic to $S$ is partitioned into
  subsets indexed by maps in $\Sigma$. The number of induced subgraphs
  that correspond to a fixed $\sigma \in \Sigma$ is given by the product
  term in Equation~(\ref{eq-ind1}).
\end{proof}

\begin{lem}
  \label{lem-gen}
  \mc{labelled} Let $\mch \subseteq \labc$ be a class of connected
  graphs that is closed under connected induced subgraphs (i.e., if a
  connected graph $F$ is in $\mch$, then every connected induced
  subgraph of $F$ is also in $\mch$). %
  Let $\iso{\mch} \coloneqq \{H_1,H_2, \dotsc\}$. Let
  $S \coloneqq \biguplus_{k=1}^m S_k$ and
  $T \coloneqq \biguplus_{k=1}^n T_k$ be two graphs such that for all
  $i$, $S_i \in \mch$ and for all $j$, $T_j \in \mch$. Let
  $\mcs \coloneqq \tup{S_1}{S_m}$. Then $\vsub{S}{T}$ and
  $\vcov{\mcs}{T}$ are functions of $\vsub{H_i}{H_j}, i,j \geq 0$.
\end{lem}

\begin{proof}
  If $\nu(S) = \nu(T)$, then $\vsub{S}{T} = 1 \text{ or } 0$ depending,
  respectively, on whether $S$ and $T$ are isomorphic or not. If
  $\nu(S) < \nu(T)$, then by Kelly's lemma \cite{kelly1942}, we have
  \[
  \vsub{S}{T} = \frac{\sum_{v\in V(T)}\vsub{S}{T-v}}{\nu(T)-\nu(S)}.
  \]
  Since each component of each induced subgraph $T-v$ is in $\mch$, an
  inductive argument on $\nu(T)$ implies that $\vsub{S}{T}$ is a
  functions of $\vsub{H_i}{H_j}, i,j \geq 0$.

  We prove by induction on $\nu(T)$ that $\vcov{\mcs}{T}$ is a functions
  of $\vsub{H_i}{H_j}, i,j \geq 0$. %
  When $\nu(T) = 1$, we have $\vsub{S}{T} = \vcov{\mcs}{T} = 1$ if
  $S \in K_1$, and $\vsub{S}{T} = \vcov{\mcs}{T} = 0$ otherwise. %
  Let the result be true when $1 \leq \nu(T) < k$. Let $\nu(T) = k$. %
  By Lemma~\ref{lem-kocay2}, we have
  \begin{equation*}
    \label{eq-gen-cov2}
    \vcov{\mcs}{T} = \prod_{i=1}^m \vsub{S_i}{T} - 
    \sum_{\substack{H \in \unlab \\
        \nu(H)\,< \,\nu(T)}}
    \vcov{\mcs}{H}\vsub{H}{T}.
  \end{equation*}
  Each factor $\vsub{S_i}{T}$ in the first term is a function of
  $\vsub{H_i}{H_j},i,j\geq 0$. %
  If there is an unlabelled graph $H$ such that $\nu(H)<\nu(T)$ that
  contributes to the summation, then $\vsub{H}{T} \neq 0$; %
  hence all connected components of $H$ are in $\mch$; %
  now by the induction hypothesis the factor $\vcov{\mcs}{H}$ is a
  function of $\vsub{H_i}{H_j},i,j\geq 0$, and $\vsub{H}{T}$ is a
  function of $\vsub{H_i}{H_j},i,j\geq 0$ as shown above. %
  This completes the induction step for $\vcov{\mcs}{T}$.
\end{proof}

In the rest of this section we show how most of the computations done
above are also possible given the abstract induced subgraph poset of a
graph. For example, to compute various invariants of a graph $G$, given
$\isp{G}$, we would need a lemma analogous to Lemma~\ref{lem-kocay2} for
tuples of elements $f_i$ of $\isp{G}$ that only implicitly specify
unlabelled graphs $F_i$. Since the induced subgraph poset does not
include empty graphs, we consider tuples $\tup{f_1}{f_k}$, where each
$f_i$ is either an empty graph $F_i = r_iK_1$ or an element of
$\isp{G}$, in which case it represents an unknown unlabelled graph
$F_i$. Our goal in the following lemma is to compute $\vcov{\mcf}{G}$
for $\mcf \coloneqq \tup{F_1}{F_k}$ given $\isp{G}$ and
$\tup{f_1}{f_k}$.

\begin{lem}{\bf (Kocay's Lemma for abstract induced subgraph posets)}
  \label{lem-kocay3}
  Let $\tup{f_1}{f_k}$ be a $k$-tuple, where each $f_i$ is either an
  empty graph $F_i = r_iK_1$ or an element of $\isp{G}$, in which case
  it represents an unknown unlabelled graph $F_i$. Let
  $\mcf \coloneqq \tup{F_1}{F_k}$. If $\nu(f_i) < \nu(G)$ for all
  $i\in\ints{1}{k}$, then $\vcov{\mcf}{G}$ can be computed from
  $\isp{G}$.
\end{lem}

\begin{proof} Let $H \vpo G$. If it is non-empty, it corresponds to an
  element $h$ of $\isp{G}$. In this case, $\nu(H) = \nu(h)$ (which is
  the rank of $h$ in $\isp{G}$, as defined earlier). For any $i$, we
  have
  \begin{equation}\label{eq-fi-h-1}
  \vsub{F_i}{H} =
  \begin{cases}
    \wv(f_i,h) & \text{ if } F_i \text{ is non-empty}, \\
    \binom{\nu(h)}{r_i} - \sum_{g\in \isp{G}\,\mid\, \nu(g) = r_i}
    \wv(g,h) & \text{ if } F_i = r_iK_1.
  \end{cases}
  \end{equation}
  On the other hand, $H$ itself could be empty, say $H = rK_1$, in which
  case
  \begin{equation}\label{eq-fi-h-2}
  \vsub{F_i}{H} =
  \begin{cases}
    0 & \text{ if } F_i \text{ is non-empty}, \\
    \binom{r}{r_i} & \text{ if } F_i = r_iK_1.
  \end{cases}
  \end{equation}
  
  If $\vcov{\mcf}{H}$ is $P$-reconstructible for each $H <_v G$, then
  $\vcov{\mcf}{G}$ is $P$-reconstructible by
  Equation~(\ref{eq-kocay2b}). We show how $\vcov{\mcf}{H}$ may be
  computed recursively on $\isp{G}$, given $\tup{f_1}{f_k}$.

  Let $r \coloneqq \max\{\nu(F_1),\dotsc,\nu(F_k)\}$. If $\nu(H) < r$,
  then $\vcov{\mcf}{H} = 0$. Therefore, we take $\nu(H) = r$ as the base
  case. There are two possibilities. If $H$ is an empty graph, then
  $\vcov{\mcf}{H} = 0$ if there are non-empty graphs in $\mcf$; if there
  are no non-empty graphs in $\mcf$, then $\vcov{\mcf}{H}$ can be
  calculated from $r_i, i = 1, \ldots, k$. If $H$ is non-empty, and $h$
  is the corresponding element of $\isp{G}$, then
  \[
  \vcov{\mcf}{H} = \prod_{i=1}^k\vsub{F_i}{H},
  \]
  where the factors on the right are computed in
  Equations~(\ref{eq-fi-h-1}) and~(\ref{eq-fi-h-2}); and since
  $\nu(H) = r = \max\{\nu(F_1),\dotsc,\nu(F_k)\}$, there are no other
  terms as in Equation~(\ref{eq-kocay2b}). Now we proceed by induction
  on $\nu(H)$ (with Equation~(\ref{eq-kocay2b})) to compute
  $\vcov{\mcf}{G}$.
\end{proof}

\begin{rem} We can now define $\vcov{\tup{f_1}{f_k}}{G}$ to mean
  $\vcov{\mcf}{G}$ whenever $\tup{f_1}{f_k}$ is a tuple whose elements
  $f_i$ are elements of the abstract induced subgraph poset $\isp{G}$ that
  correspond to unknown unlabelled graphs $F_i$, respectively, or empty
  graphs. %
  \nomenclature[4covfh]{$\vcov{\tup{f_1}{f_k}}{H}$}{similar to covers;
    $f_i$ elements of $\isp{G}$} %
\end{rem}

\begin{lem}
  \label{lem-connected}
  Connectedness of graphs is a $P$-reconstructible property.
\end{lem}

\begin{proof} We prove the claim by induction on the number of vertices
  of $G$. The base case is when $\nu(G) = 2$, in which case $G$ is
  connected if and only if $\isp{G}$ has exactly two elements. Let the
  claim be true for all graphs on at most $k$ vertices. Let
  $\nu(G) = k+1$. By induction hypothesis, we construct
  $S\coloneqq \{g_i\in \isp{G}\mid (\nu(g_i) \leq k) \wedge
  (c(G_i)=1)\}$; %
  that is, we mark elements of $\isp{G}$ of rank at most $k$ that
  correspond to connected graphs.  Now $G$ is connected if and only if
  at least two of its vertex-deleted subgraphs are connected, i.e., if
  and only if $g_M$ covers at least two distinct elements of $S$.
\end{proof}

With Lemma~\ref{lem-connected}, we assume without loss of generality
that $g_i\in \isp{G}$ are ordered so that for some $c$ the graphs
$G_i, i=\seq{1}{c}$ are connected, and the remaining graphs are
disconnected, and that $\nu(g_i) \leq \nu(g_{i+1})$ for $i \leq c-1$.

\begin{lem}
  \label{lem-cov2ref1}
  Let $\tup{s_1}{s_m} \in \set{g_1}{g_c}^m$ such that, for
  $i=\seq{1}{c}$, the element $g_i$ occurs $m_i$ times in
  $\tup{s_1}{s_m}$. It represents a graph $S = \sum_{i=1}^c
  m_iG_i$. Suppose that $\sum_{i=1}^m \nu(s_i) = \nu(G)$. Then
  \[
  \wl(S,G) = \frac{\vcov{S}{G}}{\prod_{i=1}^cm_i!},
  \]
  and hence it is $P$-reconstructible.
\end{lem}

Lemmas~\ref{lem-connected} and~\ref{lem-cov2ref1} are the basis of the
construction of the abstract bond lattice from the abstract induced
subgraph poset.

\section{Relating the abstract induced subgraph poset and the abstract
  bond lattice}
\label{sec-pg-lg}
In this section, we show that the abstract induced subgraph poset of any graph
and its abstract bond lattice can be constructed from each other except when the
graph is either a star or a disjoint union of edges. The constructions
are based on Lemma~\ref{lem-kocay3} and
Theorem~\ref{thm-thatte}-(4a). In fact, we use only the
$P$-reconstructibility of connectedness, for which we gave a very short
proof in Lemma~\ref{lem-connected}.

\begin{defn}
  \label{defn-labels}
  Let $X \coloneqq \{x_i,i \in \nn\}$ be a set of indeterminates. Let
  $\alpha\colon(\labc/{\cong})\bij X$ be a bijection. Then an unlabelled
  graph $G \coloneqq \sum_{j\in J} k_jG_j$, where $G_j, j\in J$ are
  connected unlabelled graphs, corresponds to the formal sum
  $\sum_{j\in J}k_j\alpha(G_j)$. In the notation of free abelian groups,
  let $\zz^{(X)}$ denote the set of finite formal sums of elements of
  $X$. %
  \nomenclature[4zx]{$\zz^{(X)}$}{set of finite formal sums of elements
    of $X$} %

  \medskip

  Let $f_a\colon \ispu{G}\bij \isp{G}$ be an isomorphism. %
  An {\em $X$-labelling} of $\isp{G}$ is a map
  $a\colon \isp{G}\inj \zz^{(X)}$. %
  An $X$-labelling $a$ of $\isp{G}$ %
  is {\em consistent} with $f_a$ and $\alpha$ if for all
  $G_i \in \ispu{G}\cap \labc/{\cong}$, we have
  $(af_a)(G_i) = \alpha(G_i)$; and for all $G_i\in \ispu{G}$, we have
  $(af_a)(G_i) = \sum_{j \in J}k_j\alpha(G_j)$, for some indexing set
  $J$, if and only if $G_i = \sum_{j \in J}k_jG_j$ for connected
  unlabelled graphs $G_j, j \in J$. %
  An $X$-labelling $a$ of $\isp{G}$ is {\em consistent} if there exist
  $\alpha$ and $f_a$ as described above with which $a$ is consistent.

  \medskip

  Let $f_b\colon \latu{G}\bij \lat{G}$ be an isomorphism. %
  An $X$-labelling of $\lat{G}$ is a map
  $b\colon \lat{G}\inj \zz^{(X)}$. %
  An $X$-labelling $b$ of $\lat{G}$ %
  is {\em consistent} with $f_b$ and $\alpha$ if for all
  $H_i \in \latu{G}$, $(bf_b)(H_i) = \sum_{j \in J}k_j\alpha(G_j)$, for
  some indexing set $J$, if and only if $H_i = \sum_{j \in J}k_jG_j$ for
  connected unlabelled graphs $G_j, j \in J$. %
  An $X$-labelling $b$ of $\lat{G}$ is {\em consistent} if there exist
  $\alpha$ and $f_b$ as described above with which $b$ is consistent.
\end{defn}

\begin{example} The labelling of the abstract induced subgraph poset on the right
  in Figure~\ref{fig-pg} is consistent with the map
  $\alpha(G_i) = x_i, i = 1,\ldots,5$, where the graphs $G_i$ are shown
  at the top of Figure~\ref{fig-pg}. In this example, there is a unique
  isomorphism $f_a\colon \ispu{G}\bij \isp{G}$. Similarly, a consistent
  $X$-labelling of the abstract bond lattice in Figure~\ref{fig-lc} is shown on
  the right. In this example also, the isomorphism
  $f_b\colon \latu{G}\bij \lat{G}$ is unique.
\end{example}

\subsection{From the abstract induced subgraph poset to the abstract
  bond lattice}

\begin{thm}\label{thm-pg2lg}
  The abstract bond lattice of $G$, along with a consistent $X$-labelling, can be
  constructed from its abstract induced subgraph poset.
\end{thm}

\begin{proof}
  By Lemma~\ref{lem-connected}, we construct
  $S\coloneqq \{g_i\in \isp{G}\mid c(G_i) = 1\}$; and without loss of
  generality assume that $S = \{\seq{g_1}{g_c}\}$. We define
  $\alpha(G_i) \coloneqq x_i$, for $i=\seq{1}{c}$. %
  We compute a set $B\subseteq \zz^{(X)}$ of formal sums of
  $\seq{x_1}{x_c}$ that are possible labels of $h_i \in \lat{G}$, and
  then we construct a weighted lattice on $B$ that is isomorphic to
  $\lat{G}$.

  \medskip

  \noindent {\em Computing $B$.}

  \medskip

  Let $\mcf \coloneqq (G_1^{k_1},\dotsc,G_c^{k_c})$, where $k_i \in \nn$
  are such that $\sum_i k_i\nu(G_i) = \nu(G)$. We add the formal sum
  $\sum_{i=1}^c k_ix_i$ to $B$ if and only if $\vcov{\mcf}{G} \neq 0$. %
  By Lemma~\ref{lem-kocay3}, $\vcov{\mcf}{G}$ can be computed given
  $\isp{G}$ and the tuple $(g_1^{k_1},\dotsc,g_c^{k_c})$ that
  corresponds to $\mcf$. In this manner, we compute all formal sums of
  $\seq{x_1}{x_c}$ that must be in $B$.

  \medskip

  \noindent {\em Constructing a weighted lattice structure on $B$ that
    is isomorphic to $\lat{G}$.}

  \medskip

  Let $h_a \coloneqq \sum_{i}^c a_ix_i$ and
  $h_b \coloneqq \sum_{i=1}^c b_ix_i$ be arbitrary formal sums in
  $B$. They define graphs $H_a = \sum_{i}^c a_iG_i$ and
  $H_b = \sum_{i}^c b_iG_i$. %
  By Lemma~\ref{lem-gen}, $\vcov{H_a}{H_b}$ is a function of
  $\wv(g_i,g_j), i,j\in \ints{1}{c}$ only; and
  $\wv(g_i,g_j), i,j\in \ints{1}{c}$ are known once the elements $g_i$
  in $\isp{G}$ that correspond to connected graphs are marked as in
  Lemma~\ref{lem-connected}. %
  Then by Lemma~\ref{lem-cov2ref1} we calculate
  $\wl(H_a,H_b)$. Repeating this calculation for all $h_a,h_b \in B$ we
  obtain a weighted partially ordered set with ground set $B$, that is
  isomorphic to $\latu{G}$, along with a consistent $X$-labelling on it.
\end{proof}

\subsection{From the abstract bond lattice to the abstract induced
  subgraph poset }

In this section, we construct the abstract induced subgraph poset of a graph a
from its abstract bond lattice in two steps: we first construct a consistent
$X$-labelling of the abstract bond lattice (Lemma~\ref{lem-lg2hlg}); then we
construct the abstract induced subgraph poset, along with a consistent
$X$-labelling (Lemma~\ref{lem-hlg2hpg}).
In view of Example~\ref{ex-nk2-k1n}, we exclude $K_{1,n}$ and $nK_2$ and
graphs with isolated vertices in the following result.

\begin{lem}\label{lem-lg2hlg}
  If $G$ has no isolated vertices, and $G \not \in K_{1,n}$ for any
  $n > 1$, and $G \not \in nK_2$ for any $n > 1$, then $\lat{G}$ has a
  unique consistent $X$-labelling (up to automorphisms of $\lat{G}$ and
  permutations of the indeterminates).
\end{lem}

\begin{proof}
  We construct an $X$-labelling $b\colon \lat{G}\inj \zz^{(X)}$ in two
  passes. In the beginning, we do not know the number of distinct
  connected induced subgraphs of $G$; we compute it in the first
  pass. At the same time, we construct $b(h_i)$ for each
  $h_i \in \lat{G}$, modulo the number of isolated vertices. We assign
  the minimal element $h_1$ the label $nx_1$, and assign the only
  element $h_2$ of rank 1 the label $(n-2)x_1 + x_2$, where
  $n \coloneqq \nu(G)$ is unknown. We introduce new indeterminates when
  needed. At the end of the first pass, we obtain the label
  $b(h_{\scriptscriptstyle N})$ (where $h_{\scriptscriptstyle N}$ is the
  maximal element), which, together with the assumption that $G$ has no
  isolated vertices, determines $n$. In the second pass, we calculate
  the exact labels of other elements whose labels have an $x_1$-term.

  \medskip

  \noindent {\em Assigning labels to elements $h \in \lat{G}$ of rank
    $\rho(h) = 2$.}

  \medskip

  Since $G \not \in K_{1,n}$ for any $n > 1$ and $G \not \in nK_2$ for
  any $n > 1$, both $2K_2$ and $K_{1,2}$ are subgraphs of $G$; hence
  there exist $H_3,H_4 \in \latu{G}$ of rank 2 such that
  $\{H_3,H_4\} = \{n-4)K_1+2K_2,(n-3)K_1+K_{1,2}\}$; therefore, there
  exist $h_3,h_4\in \lat{G}$ of rank 2 such that
  $\{b(h_3),b(h_4)\} = \{(n-4)x_1+2x_2, (n-3)x_1+x_3\}$. (Since we know
  that there is a connected induced subgraph from the isomorphism class
  $K_{1,2}$, we introduce a new indeterminate $x_3$ and define
  $\alpha(K_{1,2}) = x_3$). Since $\latu{K_{1,2}} \cong \latu{2K_2}$, we
  do not know $b(h_3)$ and $b(h_4)$ immediately.

  Let $\mcf \coloneqq \{S_4, C_4, K_3+K_2, K_{1,3}+K_2, P_4+K_2,
  2K_{1,2}, K_{1,2}+2K_2, K_4\}$, where $S_4$ is the graph shown in
  Figure~\ref{fig-names}. If there exist $\pi \latpo V(G)$ such that
  $\iso{G[\pi]} \in \mcf$, then $b(h_3)$ and $b(h_4)$ are uniquely
  determined. This is proved by verifying that

  \begin{enumerate}
    
  \item each graph in the above list is uniquely determined by its bond
    lattice;

  \item the abstract bond lattice of each graph in the list has only a
    trivial automorphism, which implies that the abstract bond lattice
    of each graph in the list has a unique consistent labelling up to
    choice and permutations of indeterminates.
    
  \end{enumerate}

  We demonstrate the argument for $C_4$.  Suppose that $C_4\latpo G$.
  Since $C_4$ is $\Pi$-reconstructible, there is a unique
  $h \in \lat{G}$ such that $\lat{h} \cong \latu{C_4}$. %
  Hence $b(h) = (n-4)K_1+C_4$. %
  Observe that $h$ covers exactly two elements, $h_3$ and $h_4$. Now
  $\wl((n-4)K_1+2K_2,(n-4)K_1+C_4) = 2$ and
  $\wl((n-3)K_1+K_{1,2},(n-4)K_1+C_4) = 4$. Therefore, for
  $h_i, i \in \{3,4\}$, we have $b(h_i) = (n-4)x_1+2x_2$ if
  $\wl(h_i,h) = 2$ and $b(h_i) = (n-3)x_1+x_3$ if $\wl(h_i,h) = 4$.  A
  similar argument works for all the graphs in $\mcf$.

  If there is no $F\in \mcf$ such that $F \latpo G$, then
  $\iso{G} \in \{K_4\backslash e, P_5, C_5\}$, in which case $G$ is
  $\Pi$-reconstructible; therefore, an $X$-labelling of $\lat{G}$ is
  determined up to isomorphism and choice of indeterminates. Therefore,
  in the following, we assume that there exists $F \in \mcf$ such that
  $F \latpo G$, and (without loss of generality) that
  $b(h_3) = (n-3)x_1+x_3$ and $b(h_4) = (n-4)x_1+2x_2$.

  \medskip

  \noindent {\em Assigning labels to elements $h \in \lat{G}$ of rank
    $\rho(h) > 2$ in the case when the corresponding graph $H$ has at
    most one non-trivial component.}

  \medskip

  For each connected graph $G_i\in \ispu{G}$, there is a distinct graph
  $H_i \in \latu{G}$ such that $H_i = m_iK_1 + G_i$. Hence we define
  $\Gamma \coloneqq \{h_i \in \lat{G} \mid H_i = m_iK_1 + G_i \text{ for
    some connected graph } G_i\}$,
  and $\Gamma_r \coloneqq \{h\in \Gamma\mid \rho(h) \leq r\}$, and
  construct $\Gamma$ by constructing $\Gamma_r$ recursively.

  We have $\Gamma_2 = \{h_1,h_2, h_3\}$. Suppose that we have
  constructed $\Gamma_k$ for all $k\in\{2,\dotsc,r\}$. Let
  $h_i\in \lat{G}$ such that $\rho(h_i) = r+1$. We use the fact that a
  graph on two or more vertices is connected if and only if at least 2
  of its vertex-deleted subgraphs are connected. If $h_i$ is in
  $\Gamma_{r+1}$, then $h_i$ must cover at least two distinct elements
  $h_j, h_k$ in $\Gamma_{r}$. But the converse is not true. If $h_i$
  covers two distinct elements $h_j, h_k$ in $\Gamma_{r}$, then either
  $H_i = aK_1+K_2+F$ for some $a\geq 0$ and a connected unlabelled
  graphs $F$, or $H_i= aK_1+F$ for some $a\geq 0$ and a connected
  unlabelled graph $F$. We want to add $h_i$ to $\Gamma_{r+1}$ only in
  the latter case. The necessary and sufficient condition for the former
  case is that there exists an element $h_\ell$ covered by $h_i$ such
  that $\epsilon(H_\ell) = \epsilon(H_i)-1$ and
  $\wl((n-4)K_1+2K_2,H_i) = \wl((n-4)K_1+2K_2,H_\ell) +
  \epsilon(H_\ell)$.
  These conditions are recognised from $\lat{G}$ and the induction
  hypothesis since $\epsilon(H_\ell) = \wl(h_2,h_\ell)$. If such an
  element $h_\ell$ does not exist, then we add $h_i$ to $\Gamma_{r+1}$
  and define $b(h_i) \coloneqq (n-r-2)x_1+x_i$, where we have introduced
  a new indeterminate $x_i$ for the (unknown) graph $F$ such that
  $H_i= (n-r-2)K_1+F$.

  Without loss of generality, we assume that
  $b(h_i) = (n-\rho(h_i)-1)x_1 + x_i$, for $i=\seq{1}{c}$, and that
  $\seq{h_1}{h_c} \in \Gamma$ are ordered so that
  $\rho(h_i) \leq \rho(h_{i+1})$ for all $i \in \ints{1}{c-1}$.

  \medskip

  \noindent {\em Assigning labels to elements $h_i\in \lat{G}, i > c$
    (i.e., to $h_i \in \lat{G}\backslash \Gamma$) in terms of
    $\seq{x_1}{x_c}$.}

  \medskip

  Recall the notation $c(G_i,G_j)$ for the number of components of $G_j$
  that are isomorphic to $G_i$. We have
  \begin{align}
    \label{eq-lg2hlg}
    \vsub{G_j}{H_i} 
    & = c(G_j,H_i) + \sum_{r = j+1}^c \vsub{G_j}{G_r}c(G_r,H_i) 
      \text{ for all } j \leq c \text{ and for all } i > c\nonumber \\
    \therefore \wl(h_j, h_i) 
    &= c(G_j,H_i) + \sum_{r=j+1}^c \wl(h_j,h_r)c(G_r,H_i)  
      \text{ for all } j \leq c \text{ and for all } i > c.
  \end{align}

  Now $b(h_i) = \sum_{j=1}^c x_jc(G_j,H_i)$ is obtained for each
  $h_i, i > c$ by solving the system of equations~(\ref{eq-lg2hlg}) for
  $c(G_j,H_i)$ for all $j \leq c$. Indeed, for all $i > c$,
  $c(G_c,H_i) = \wl(h_c, h_i)$, and for all $j < c$, for all $i > c$,
  $c(G_j,H_i)$ is expressed in terms of $c(G_r,H_i), r > j$.

  \medskip

  \noindent {\em Inferring $n$ from the label of the maximum element
    $h_{\scriptscriptstyle N}$, and fixing labels that contain $n$.}

  \medskip

  Once $b(h_{\scriptscriptstyle N}) = \sum_{i=2}^ck_ix_i$ has been
  calculated, we calculate
  $n = \nu(G) = \rho(G) + c(G) = \rho(h_{\scriptscriptstyle N}) +
  \sum_{i=2}^ck_i$;
  then for all $h\in \lat{G}$, the $x_1$ term in $b(h)$ is determined
  exactly.
\end{proof}

\begin{lem}\label{lem-hlg2hpg}
  The abstract induced subgraph poset of $G$ along with a consistent
  $X$-labelling of $G$ can be constructed from its consistently
  $X$-labelled abstract bond lattice.
\end{lem}

\begin{proof} Let $b \colon \lat{G}\inj \zz^{(X)}$ be a consistent
  $X$-labelling of $\lat{G}$. Without loss of generality, we assume that
  $b(h_i)= (n-n_i)x_1+x_i$, for $i=\seq{1}{c}$, and $b(h_i), i > c$ are
  other distinct formal sums of $x_i, i=\seq{1}{c}$, and that $n-n_i$
  are in non-increasing order. We first construct a set
  $A \subseteq \zz^{(X)}$ of possible formal sums that can be in an
  $X$-labelling $a$ of $\isp{G}$. Then we make $A$ a weighted partially
  ordered set that is isomorphic to $\isp{G}$.
  
  We have $\set{x_1}{x_c} \subseteq A$; and
  $b(h_{\scriptscriptstyle N}) \in A$ since
  $a(g_{\scriptscriptstyle M}) = b(h_{\scriptscriptstyle N})$. %
  If $b(h_1)= nx_1$, then $\nu(G) = n$; hence $\nu(g_i)\coloneqq n_i$ is
  determined by $b(h_i)$ for all $i\in\ints{1}{c}$.
  Moreover, for all $i,j$ such that $1 \leq i \leq j \leq c$,
  \begin{equation}
    \label{eq-hi-hj}
    w_v(x_i,x_j) \coloneqq \vsub{G_i}{G_j} =
    \begin{cases}
      n_j & \text{ if } i = 1, \\
      \wl(h_i,h_j) & \text{ otherwise.} \\
    \end{cases}
  \end{equation}

  Let $\mch \coloneqq \cup_{i=1}^c G_i$ (where $G_i$ are isomorphism
  classes). Lemma~\ref{lem-gen} may be applied to $\mch$ since it is
  closed under connected induced subgraphs, and $\vsub{G_i}{G_j}$ are
  calculated above. Let $F = \sum_i k_iG_i$ for some $k_i \in \nn$. By
  Lemma~\ref{lem-gen}, $\vsub{F}{G}$ is a function of $w_v(x_i,x_j)$,
  $a(g_{\scriptscriptstyle M})$, and $k_i, i=\seq{1}{c}$. So we can
  determine if any given linear combination $\sum_{i=1}^ck_ix_i$ is in
  $A$. Thus we list all linear combinations $\sum_{i=1}^ck_ix_i$ that
  must be in $A$. Then again by Lemma~\ref{lem-gen} we calculate
  $\vsub{G_j}{G_k}$ for any two graphs $G_j,G_k \in \ispu{G}$
  represented, respectively, by linear combinations
  $\sum_{i=1}^c k_ix_i$ and $\sum_{i=1}^c l_ix_i$ in $\isp{G}$. This
  completes the construction of a consistently $X$-labelled $\isp{G}$.
\end{proof}

\begin{thm}\label{thm-lg2pg}
  If $G$ has no isolated vertices, and $G \not \in K_{1,n}$ for any
  $n > 1$, and $G \not \in nK_2$ for any $n > 1$, then $\isp{G}$ can be
  constructed from $\lat{G}$.
\end{thm}

\begin{proof}
  Either $G$ itself is $\Pi$-reconstructible (as in the case of
  $K_4\backslash e$, $P_5$ and $C_5$), or the claim follows from
  Lemmas~\ref{lem-lg2hlg} and~\ref{lem-hlg2hpg}.
\end{proof}

\begin{cor}\label{cor-ulam}
  \text{}
  \begin{enumerate} 
  \item Ulam's conjecture is true if and only if every graph
    $G \in \lab$ such that $G$ has no isolated vertices, and
    $G \not \in K_{1,n}$ for any $n > 1$, and $G \not \in nK_2$ for any
    $n > 1$, can be constructed up to isomorphism from its abstract bond
    lattice $\lat{G}$.
  
  \item\label{cor-ulam-trees} Every tree or forest $T$ such that
    $T \not \in K_{1,n}$ for any $n > 1$, $T \not \in nK_2$ for any
    $n > 1$, and $T$ has no isolated vertices, can be constructed up to
    isomorphism from its abstract bond lattice $\lat{T}$.
  \end{enumerate}
\end{cor}

\begin{proof}
  The first part follows from Theorem~\ref{thm-lg2pg} and
  Theorem~\ref{thm-thatte}-(\ref{thm-thatte-ulam}). The second part of
  the theorem now follows from Theorem~\ref{thm-lg2pg} and
  Theorem~\ref{thm-thatte}-(\ref{thm-thatte-trees}).
\end{proof}

\section{Colouring polynomials of graphs}
\label{sec-polynomials}

We define two symmetric polynomial invariants of graphs, namely, the
{\em chromatic symmetric function} and a stronger invariant called the
{\em symmetric Tutte polynomial}, which were introduced by
\citet{stanley1995,stanley1998}.

A {\em vertex colouring} of a graph $G$ is a map
$\kappa: V(G)\to \zz^+$. %
\nomenclature[5k]{$\kappa$}{vertex colouring function on a graph} %
Let $\kappa$ be a vertex colouring of $G$, and let $U \subseteq
V(G)$.
We say that $U$ (or a subgraph on $U$) is {\em monochromatic} if
$\kappa$ is constant on $U$. We say that $\kappa$ is {\em proper} if
there are no monochromatic edges. %

Let $x_1,x_2,\dotsc$ and $t$ be commuting indeterminates. %
\nomenclature[4x3]{$x_1,x_2,\dotsc$ and $t$}{commuting indeterminates} %
We denote by $1^n$ the assignment $x_1=x_2=\cdots=x_n=1$ and
$x_{n+1}=x_{n+2}=\cdots=0$. We write $x_1,x_2,\dotsc$ collectively as
just $x$. %
\nomenclature[21]{$1^n$}{assignment $x_i=1$ for $i \leq n$, and $x_i=0$
  for $i> n$}

\begin{defn}\label{def-xg} The
  {\em chromatic symmetric function} $X_G(x)$ of $G$ is defined by
  \[
  X_G(x) \coloneqq \sum_\kappa \prod_{v \in V(G)} x_{\kappa(v)},
  \]
  where the summation is over all proper colourings $\kappa$ of $G$. %
  \nomenclature[4x2gx]{$X_G(x)$}{chromatic symmetric function of a graph
    $G$}
\end{defn}

\begin{defn}\label{def-xtg} The {\em symmetric Tutte polynomial} of $G$
  is defined by
  \[
  X_G(x;t) \coloneqq \sum_\kappa (1+t)^{m(\kappa)} %
  \prod_{v \in V(G)} x_{\kappa(v)},
  \]
  where the summation is over all vertex colourings $\kappa$ of $G$, and
  $m(\kappa)$ is the number of monochromatic edges of $\kappa$.
  \nomenclature[4mk]{$m(\kappa)$}{number of monochromatic edges in a
    vertex colouring $\kappa$} %
  \nomenclature[4x2gxt]{$X_G(x;t)$}{symmetric Tutte polynomial of a
    graph $G$} %
\end{defn}

\citet{stanley1995} and \citet{chow1995} have noted that
$X_G(x;-1) = X_G(x)$; and $X_G(1^n)$ and $X_G(1^n;t)$ are equivalent to
the chromatic polynomial and the Tutte polynomial of $G$, respectively,
where $n \coloneqq \nu(G)$. Thus the symmetric Tutte polynomial
$X_G(x;t)$ specialises to the other three invariants mentioned above.

\begin{thm}[\citet{stanley1995,stanley1998} and
  \citet{chow1995}]\label{thm-stanley-xgt} We have
  \[
  X_G(x;t) = \sum_{S\subseteq E(G)} t^{|S|} \sum_{\kappa} %
  \prod_{v \in V(G)} x_{\kappa(v)},
  \]
  and
  \[
  X_G(x) = X_G(x,-1) = \sum_{S\subseteq E(G)} (-1)^{|S|} \sum_{\kappa} %
  \prod_{v \in V(G)} x_{\kappa(v)},
  \]
  where the inner summation in each equation is over all vertex
  colourings $\kappa$ that are monochromatic on the connected components
  of the spanning subgraph $G_S$ of $G$ with edge set $S$.
\end{thm}

\subsection{Reconstructing the symmetric Tutte polynomial}
\label{sec-xtg}

The reconstructibility of $X_G(x)$ and $X_G(x;t)$ is not immediately
obvious from their expansions given in Theorem~\ref{thm-stanley-xgt}. In
the following, we prove, using Kocay's lemma for abstract induced
subgraph poset (Lemma~\ref{lem-kocay3}), that $X_G(x;t)$ is
reconstructible (in fact $P$-reconstructible), which implies that the
chromatic symmetric function is reconstructible as well.

\begin{thm}\label{thm-xgt}
  The symmetric Tutte polynomial is $P$-reconstructible.
\end{thm}

\begin{proof}
  Each colouring $\kappa: V(G) \to \zz^+$ with $k$ colours defines two
  $k$-tuples: a tuple $\tup{m_1}{m_k}$ of positive integers, where
  $m_1 < \cdots < m_k$ are the $k$ colours, and a tuple $\tup{F_1}{F_k}$
  of vertex-disjoint subgraphs $F_i$ of $G$, which are induced by the
  colour classes $m_i$, respectively; hence $\sum_i\nu(F_i) =
  \nu(G)$. Therefore, by the definition of $X_G(x;t)$, we have
  \[
  X_G(x;t) = \sum_{k\in \zz^+}\;
  \sum_{\substack{\tup{m_1}{m_k}\mid\\m_1<\cdots <m_k}}\;
  \sum_{\tup{F_1}{F_k}} (1+t)^{\sum_i \epsilon(F_i)}
  \prod_{i}^kx_{m_i}^{\nu(F_i)}.
  \]

  We say that tuples $\tup{F_1}{F_k}$ and $\tup{F^\prime_1}{F^\prime_k}$
  are {\em equivalent} if there is a bijection
  $f\colon \ints{1}{k}\to\ints{1}{k}$ such that
  $F_i \cong F^\prime_{f(i)}$ for all $i$; else they are {\em
    inequivalent}. Similarly, we say that tuples $\tup{f_1}{f_k}$ and
  $\tup{f_1^\prime}{f_k^\prime}$, where $f_i$ and $f_i^\prime$ are
  elements of $\isp{G}$ or empty graphs, are {\em equivalent} if there
  is a bijection $g\colon \ints{1}{k}\to\ints{1}{k}$ such that
  $f_i = f^\prime_{g(i)}$ for all $i$; else they are {\em inequivalent}.

  We write the inner summation over mutually inequivalent tuples
  $\tup{F_1}{F_k}$. There are $\vcov{\tup{F_1}{F_k}}{G}$ tuples in the
  equivalence class of any given tuple $\tup{F_1}{F_k}$. Hence
  \[
  X_G(x;t) = \sum_{k\in \zz^+}\;
  \sum_{\substack{\tup{m_1}{m_k}\mid\\m_1<\cdots <m_k}}\;
  \sum_{\substack{\text{inequivalent}\\ \tup{F_1}{F_k}}}
  \vcov{\tup{F_1}{F_k}}{G} (1+t)^{\sum_i \epsilon(F_i)}
  \prod_{i}^kx_{m_i}^{\nu(F_i)}.
  \]

  Finally, by Lemma~\ref{lem-kocay3}, we compute the inner summation
  over mutually inequivalent tuples $\tup{f_1}{f_k}$, where each $f_i$
  is an element of $\isp{G}$ or is an empty graph, and
  $\sum_i\nu(f_i) = \nu(G)$. Hence
  \[
  X_G(x;t) = \sum_{k\in \zz^+}\;
  \sum_{\substack{\tup{m_1}{m_k}\mid\\m_1<\cdots <m_k}}\;
  \sum_{\substack{\text{inequivalent}\\ \tup{f_1}{f_k}}}
  \vcov{\tup{f_1}{f_k}}{G} (1+t)^{\sum_i \epsilon(f_i)}
  \prod_{i}^kx_{m_i}^{\nu(f_i)}.
  \]
  This implies the $P$-reconstructibility of $X_G(x;t)$.
\end{proof}

In fact, $X_G(x;t)$ may be generalised as follows. Define a map
$\kappa: V(G)\to 2^{\zz^+}$ to be a {\em set-colouring} of $G$ when each
colour set $\kappa(v)$ is finite. The invariant $Y_G(x;y)$ over
indeterminates $x_1,x_2,\dotsc$ and $y_i, y_2,\dotsc$ defined below
generalises the symmetric Tutte polynomial in a natural way:
\[
Y_G(x;y) \coloneqq \sum_{\kappa} \prod_{\{u,v\} \in E(G)} \; \prod_{i\in
  \kappa(u)\cap \kappa(v)} (1+y_i) \quad \prod_{w\in V(G)} \;\prod_{j\in
  \kappa(w)}x_j.
\]
We call it the {\em set-colouring symmetric Tutte polynomial} of $G$. We
do not know if the invariant $Y_G(x;y)$ is strictly stronger than
$X_G(x;t)$ (i.e., if there are graphs distinguished by this stronger
invariant but not by their symmetric Tutte polynomial). In this paper we
do not study it beyond observing (without proof) that $Y_G(x;y)$ is
$P$-reconstructible; the proof is similar to that of
Theorem~\ref{thm-xgt}, except that we do not require that
$\sum_{i}\nu(f_i) = \nu(G)$ since the colour classes of a colouring do
not partition $V(G)$.

For other known proofs of the reconstructibility of the various
colouring polynomials, we refer to the induced subgraph expansion of the
chromatic polynomial in \cite{biggs1993}, the reconstructibility of the
chromatic symmetric function in \cite{chow1995}, and the
reconstructibility of the chromatic and Tutte polynomials in
\cite{bondy1991}. As far as the author is aware, the reconstructibility
of the symmetric Tutte polynomial presented above is a new result.

\subsection{Reconstructing the chromatic symmetric function}
\label{sec-xg}
The reconstructibility of the chromatic symmetric function follows from
the reconstructibility of the symmetric Tutte polynomial. But here we
given another proof based on its expansion given by \citet[Theorem
2.6]{stanley1995} (which is Theorem~\ref{thm-xg-stanley} in this
section) and our Theorem~\ref{thm-pg2lg}.

Let $\zeta$ and $\mu$ be, respectively, the zeta function and the
M\"obius function of $\latl{G}$ (which is the (unfolded) connected
partition lattice of $G$ defined in Section~\ref{sec-intro-lattice}). %
\nomenclature[5f]{$\zeta(.,.)$}{zeta function of a poset} %
\nomenclature[5m]{$\mu(.,.)$}{M\"obius function of a poset} %

By definition, for any two connected partitions $\pi$ and $\tau$ of
$V(G)$, we have
\begin{equation}
  \label{eq-mu-zeta}
  \sum_{\sigma} \mu(\pi, \sigma)\zeta(\sigma,\tau) =
  \begin{cases}
    1 & \text{ if } \tau = \pi \\
    0 & \text{ otherwise}.
  \end{cases}
\end{equation}
Moreover, if $\sigma_1$ and $\sigma_2$ are any two connected partitions
of $V(G)$ such that $G[\sigma_1] \cong G[\sigma_2]$, then
$\mu(\hat{0},\sigma_1) = \mu(\hat{0},\sigma_2)$. Therefore, a function
$\psi\colon \lat{G} \to \zz$ given by
$\psi(h_i) \coloneqq \mu(\hat{0},\sigma)$, where $\sigma$ is any
connected partition of $V(G)$ such that $G[\sigma]\cong H_i$, is well
defined. %
\nomenclature[5z2]{$\psi(.)$}{a function related to the M\"obius function
  of a bond lattice} %
Therefore, when $\pi = \hat{0}$, the system of
equations~(\ref{eq-mu-zeta}) may be re-written as
\begin{equation}
  \label{eq-psi-w}
  \sum_{h_i} \psi(h_i)w_\sigma(h_i,h_j) =
  \begin{cases}
    1 & \text{ if } h_j=h_1 \\
    0 & \text{ otherwise},
  \end{cases}
\end{equation}
for all $h_j \in \lat{G}$.

\begin{lem}\label{lem-psi} 
  The abstract bond lattice $\lat{G}$ uniquely determines
  $\psi(h_i), h_i \in \lat{G}$.
\end{lem}
\begin{proof}
  The weight function $w_\sigma$ is invertible in the incidence algebra
  of $\lat{G}$. Therefore, the system of equations~(\ref{eq-psi-w}) has
  a unique solution for $\psi(h_i), h_i \in \lat{G}$.
\end{proof}

For an integer partition $\lambda \coloneqq \tup{\lambda_1}{\lambda_l}$,
let $p_{\lambda}$ denote its {\em power sum symmetric function} given by
$p_{\lambda} = \prod_{k=1}^l\sum_ix_i^{\lambda_k}$,
see \citet{stanley-v1}. %
\nomenclature[4pl]{$p_{\lambda}$}{power sum symmetric function of
  integer partition $\lambda$}
\begin{thm}[\citet{stanley1995}]\label{thm-xg-stanley} We have
  \begin{equation}
    \label{eq-xg-1}
    X_G(x) = \sum_{\pi \in \latl{G}}\mu(\hat{0},\pi)p_{\lambda(\pi)}.
  \end{equation}
\end{thm}

\begin{cor}\label{cor-xg3}
  If $G$ has no isolated vertices, and $G \not \in K_{1,n}$ for any $n >
  1$, and $G \not \in nK_2$ for any $n > 1$, then $X_G(x)$ can be
  computed from $\lat{G}$.
\end{cor}

\begin{proof}
  Equation~(\ref{eq-xg-1}) may be re-written as
  \begin{equation}\label{eq-xg-2}
    X_G(x) = \sum_{h_i\in\lat{G}}\psi(h_i)
    w_\sigma(h_i,h_{\scriptscriptstyle N})p_{\lambda(h_i)},
  \end{equation}
  where $\lambda(h_i) \coloneqq \lambda(H_i)$. For each $h_i$ in
  $\lat{G}$, $\lambda(H_i)$ is determined by Theorem~\ref{thm-pg2lg}
  and $\psi(h_i)$ is determined by Lemma~\ref{lem-psi}.
\end{proof}

\subsection{Colouring polynomials of trees}
\label{sec-trees}

Let $T$ be a tree. Let
$\lat{T}\coloneqq \{h_1,\dotsc,h_{\scriptscriptstyle N}\}$, where $h_i$
are enumerated so that the number of components $c(H_i)$ are
non-increasing. Let $\mu$ be the M\"obius function of $\latl{T}$, and
let $\psi\colon \lat{T} \to \zz$ be as defined in Section~\ref{sec-xg}.

\begin{lem}
  \label{lem-mu-tree} 
  For all $h_i \text{ in } \lat{T}$,
  \begin{equation}
    \psi(h_i) = (-1)^{\epsilon(H_i)} = (-1)^{\nu(T)-c(H_i)}. 
  \end{equation}
\end{lem}

\begin{proof}
  The connected partition lattice of $T$ is isomorphic to the power set
  lattice of $E(T)$. The M\"obius function of the power set lattice of
  $E(T)$ is given by $\mu(E_1,E_2) = (-1)^{|E_2|-|E_1|}$, for
  $E_1,E_2 \subseteq E(T)$. Now the result follows from the definition
  of $\psi$.
\end{proof}

Given a graph $G$ and a partition $\lambda$ of $\nu(G)$, let
$k(\lambda,G) \coloneqq |\{\pi \models_c V(G)\mid \lambda(\pi) =
\lambda\}|$,
and let
$\parts(G) \coloneqq \{(\lambda, k(\lambda,G)) \mid \lambda\models
\nu(G)\}$.
\nomenclature[4k3lg]{$k(\lambda,G)$}{number of connected partitions of
  $G$ of type $\lambda$} %
\begin{lem}\label{lem-xg2partitions}
  The invariants $\parts(T)$ and $X_T(x)$ are equivalent, i.e., they can
  be computed from each other.
\end{lem}

\begin{proof}
  By Equation~(\ref{eq-xg-2}) and Lemma~\ref{lem-mu-tree}, we have
  \begin{equation}
    \label{eq-klambda1}
    X_T(x)
    = \sum_{\lambda \,\vdash \,\nu(T)}(-1)^{\nu(T)-\ell(\lambda)} k(\lambda,T)\,p_{\lambda},
  \end{equation}
  which shows that $X_T(x)$ can be computed from $\parts(T)$.

  Given $X_T(x)$, equation~(\ref{eq-klambda1}) can be solved for
  $k(\lambda,T)$ as follows. Let
  $\lambda \coloneqq \tup{\lambda_1}{\lambda_l}$ be a partition of
  $\nu(T)$, and let $b_T(\lambda)$ be the coefficient of
  $\prod_{i=1}^lx_i^{\lambda_i}$ in $X_T(x)$. The numbers $b_T(\lambda)$
  and $k(\lambda,T)$ satisfy the equation
  \begin{equation}
    \label{eq-klambda2}
    b_T(\lambda) = 
    \sum_{\lambda^{\prime}\models \lambda}
    (-1)^{\nu(T)-l(\lambda^{\prime})}k(\lambda^{\prime},T)
    a(\lambda^{\prime},\lambda),
  \end{equation}
  where $a(\lambda^{\prime},\lambda)$ is the coefficient of
  $\prod_{i=1}^lx_i^{\lambda_i}$ in $p_{\lambda^{\prime}}$. Moreover,
  $a(\lambda^{\prime},\lambda) = 0$ if $\lambda^\prime$ is not a
  refinement of $\lambda$, and $a(\lambda^{\prime},\lambda) = 1$ if
  $\lambda^\prime = \lambda$.  The system of
  equations~(\ref{eq-klambda2}) can be recursively solved for
  $k(\lambda,T)$ (the initial condition being $k(\lambda,T)=1$ when
  $\lambda$ is the finest partition of $\nu(T)$). Thus we can construct
  $\parts(T)$ given $X_T(x)$.
\end{proof}

\begin{figure}[ht]
  \begin{center}
    \includegraphics{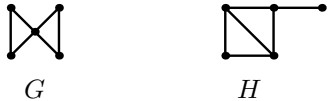}
  \end{center}
  \caption[]{Non-isomorphic graphs with identical chromatic symmetric
    function.}
  \label{fig-stan}
\end{figure}

\begin{rem}
  The above result is not valid for graphs in general; e.g., the graphs
  $G$ and $H$ shown in Figure~\ref{fig-stan} have the same chromatic
  symmetric function (see \citet{stanley1995}), but we can verify that
  $\parts(G)$ and $\parts(H)$ are different.
\end{rem}

\begin{lem}
  \label{lem-recognition} Whether $G$ is a tree or not can be recognised
  from $X_G(x)$ and from $\parts(G)$.
\end{lem}

\begin{proof}
  We obtain $\nu(G)$ and $\epsilon(G)$ by observing that
  $\nu(G)$ is the degree of each monomial in $X_G(x)$; it is the length
  of the finest partition $(1,1,\dotsc,1)$ that appears in a pair in
  $\parts(G)$;
  the coefficient of $x_1^2x_2x_3\dotsm x_{\nu(G)-1}$ in $X_G(x)$ is
  $\binom{\nu(G)}{2}-\epsilon(G)$; and
  for $\lambda \coloneqq (2,1,1,\dotsc,1)$, we have
  $k(\lambda,G)=\epsilon(G)$.

  Now $G$ is a tree if and only it $\epsilon(G) = \nu(G)-1$ and it has
  no cycles; and $G$ is acyclic if and only if the number of acyclic
  orientations of $G$ is $2^{\epsilon(G)}$. \citet{stanley1973} has
  shown that the number of acyclic orientations of a graph $G$ can be
  calculated from its chromatic polynomial $\chi_G(x)$, which is given
  by $X_G(1^n)$.

  The number of components in $G$ is $\min(\ell(\lambda))$ over
  $(\lambda, k(\lambda,G))$ in $\parts(G)$. Now $G$ is a tree if and
  only if $\epsilon(G) = \nu(G)-1$ and $G$ has a single component.
\end{proof}

\begin{notation} Let $\mathbf{v} \coloneqq \tup{v_1}{v_r}$ and
  $\mathbf{e} \coloneqq \tup{e_1}{e_r}$ and
  $\lambda\coloneqq \tup{\lambda_1}{\lambda_s}$ be integer vectors. %
  Let
  \[
  \phi(\mathbf{v},\mathbf{e},\lambda) \coloneqq \left|\left\{f\colon
      \ints{1}{s}\sur\ints{1}{r} \,\middle|\, \left(\forall j \in
        \ints{1}{r}\right) \left(\sum_{i \in f^{-1}(j)}\lambda_i = v_j,
        \sum_{i \in f^{-1}(j)}(\lambda_i-1) = e_j\right)\right\}\right|.
  \]
  \nomenclature[5z1]{$\phi(\mathbf{v},\mathbf{e},\lambda)$}{certain
    function of integer vectors $\mathbf{v},\mathbf{e}$ and an integer
    partition $\lambda$} %
  Let $\theta(\mathbf{v},\mathbf{e};T)$ be the number of ordered
  partitions $\tup{V_1}{V_r}$ of $V(T)$ such that $|V_i| = v_i$, and
  $\epsilon(T[V_i]) = e_i$; we call such partitions
  $(\mathbf{v},\mathbf{e})$-partitions of $V(T)$. %
  \nomenclature[5theta]{$\theta(\mathbf{v},\mathbf{e};T)$}{number of
    certain partitions defined for a tree} %
  For vectors $\mathbf{v}$ and $\mathbf{e}$ of unequal length, we define
  $\phi(\mathbf{v},\mathbf{e},\lambda)\coloneqq 0$ and
  $\theta(\mathbf{v},\mathbf{e};T) \coloneqq 0$. %
  Given integer vectors $\mathbf{e} \coloneqq \tup{e_1}{e_r}$ and
  $\mathbf{f}\coloneqq \tup{f_1}{f_s}$, we define
  $\mathbf{e}\leq \mathbf{f}$ if and only if $r=s$ and $e_i \leq f_i$
  for all $i$; the relation $\leq$ makes the set of integer vectors a
  locally finite partially ordered set.
\end{notation}

\begin{lem}
  \label{lem-theta}
  For all integer vectors $\mathbf{v}$ and $\mathbf{e}$, the parameter
  $\theta(\mathbf{v},\mathbf{e};T)$ is determined by $X_T(x)$.
\end{lem}

\begin{proof}
  For all vectors $\mathbf{v}$ and $\mathbf{e}$, we have
  \begin{equation}
    \label{eq-theta}
    \sum_{\mathbf{f} \mid \mathbf{e}\leq \mathbf{f}}
    \theta(\mathbf{v},\mathbf{f};T)\prod_{i=1}^r \binom{f_i}{e_i}
    = \sum_{\pi\,\models_c\, V(T)}\phi(\mathbf{v},\mathbf{e},\lambda(\pi))
    = \sum_{\lambda\,\models\, \nu(T)}\phi(\mathbf{v},\mathbf{e},\lambda)k(\lambda,T).
  \end{equation}
  Equation~(\ref{eq-theta}) is trivially true if vectors $\mathbf{v}$
  and $\mathbf{e}$ have unequal lengths, or $\mathbf{v}$ has some
  non-positive entries, or $\mathbf{e}$ has some negative
  entries. Otherwise, we prove the equation by double counting as
  follows. %
  Let vectors vectors $\mathbf{v}$ and $\mathbf{e}$, both of length $r$,
  be fixed. We define a matrix $A \coloneqq A_{\mathbf{v},\mathbf{e}}$
  with rows indexed by connected partitions of $V(T)$, and columns
  indexed by $(\mathbf{v},\mathbf{f})$-partitions of $V(T)$ such that
  $\mathbf{f} \geq \mathbf{e}$. For a connected partition
  $\pi \coloneqq \set{X_1}{X_s}$ of $V(T)$ and a
  $(\mathbf{v},\mathbf{f})$-partition $\sigma \coloneqq \tup{Y_1}{Y_r}$
  of $V(T)$, an entry $A(\pi,\sigma)$ is 1 if $\pi$ refines $\sigma$,
  and for all $j \in \ints{1}{r}$,
  $\sum_{X_i\mid X_i \subseteq Y_j} |X_i| = v_j$ and
  $\sum_{X_i\mid X_i \subseteq Y_j} (|X_i|-1) = e_j$; and
  $A(\pi,\sigma)$ is 0 otherwise. %
  For a fixed vector $\mathbf{f}$, the number of 1s in a column of $A$
  indexed by a $(\mathbf{v},\mathbf{f})$-partition is
  $\prod_{i=1}^r \binom{f_i}{e_i}$, and the number of
  $(\mathbf{v},\mathbf{f})$-partitions is
  $\theta(\mathbf{v},\mathbf{f};T)$. %
  The number of 1s in a row indexed by a connected partition
  $\pi\coloneqq \set{X_1}{X_s}$ is
  $\phi(\mathbf{v},\mathbf{e},\lambda(\pi))$. %
  Now counting the number of 1s in $A$ by columns and rows gives the
  first equality.  The second equality follows from the definition of
  $k(\lambda,T)$.
  
  By Lemma~\ref{lem-xg2partitions}, $\parts(T)$ can be constructed from
  $X_T(x)$; therefore, the \rhs of Equation~(\ref{eq-theta}) is
  known. Now Equation~(\ref{eq-theta}) is solved for
  $\theta(\mathbf{v},\mathbf{e};T)$ for all $\mathbf{v}$ and
  $\mathbf{e}$ by M\"obius inversion on the partially ordered set of
  integer vectors.
\end{proof}

\begin{thm}
  \label{lem-xgt-tree} For trees, the chromatic symmetric function and
  the symmetric Tutte polynomial are equivalent.
\end{thm}

\begin{proof}
  For all graphs, $X_G(x)$ is a specialisation of $X_G(x;t)$. We show
  how $X_G(x;t)$ is determined by $X_G(x)$ when $G$ is a tree. We have
  \begin{equation}
    \label{eq-xgt}
    X_T(x;t) = 
    \sum_{k\in \zz^+}\quad \sum_{\mathbf{v_k},\mathbf{e_k},\mathbf{m_k}}
    \theta(\mathbf{v_k},\mathbf{e_k};T)
    (1+t)^{\sum_i e_i}\prod_{i=1}^{k}x_{m_i}^{v_i},
  \end{equation}
  where
  $\mathbf{v_k}\coloneqq \tup{v_1}{v_k} \in \zz^k,\mathbf{e_k}\coloneqq
  \tup{e_1}{e_k} \in \zz^k$,
  $\mathbf{m_k}\coloneqq \tup{m_1}{m_k} \in \zz^k$ are such that
  $v_i > 0, e_i \geq 0$ for all $i$, and $0 < m_1 < \cdots < m_k$. By
  Lemma~\ref{lem-theta}, $\theta(\mathbf{v},\mathbf{e};T)$ is determined
  by $X_T(x)$ for all $\mathbf{v}$ and $\mathbf{e}$. Thus, for trees,
  the two invariants are equivalent.
\end{proof}

\begin{rem}
  Theorem~\ref{lem-xgt-tree} implies that Stanley's question (whether
  the chromatic symmetric function distinguishes trees) and the question
  of Noble and Welsh (whether their weighted chromatic function with
  unit weights, which is equivalent to the symmetric Tutte polynomial,
  distinguishes trees) are equivalent.
\end{rem}

Lemma~\ref{lem-theta} may be a useful tool to study Stanley's question.
We illustrate one simple application of the lemma. Let the {\em degree
  of a subtree} $F$ of a tree $T$ be the number of edges in $T$ with one
end in $F$ and one end outside $F$.

\begin{cor}
  \label{lem-deg} The number of subtrees of $T$ with a given number of
  vertices and a given degree is determined by $X_T(x)$; in particular,
  the degree sequence of $T$ is determined by $X_T(x)$.
\end{cor}

\begin{proof}
  Let $\mathbf{v} \coloneqq (k, \nu(T)-k)$ and
  $\mathbf{e} \coloneqq (k-1, \nu(T)-k -d)$.  Now the number
  $\theta(\mathbf{v}, \mathbf{e};T)$ counts the number of subtrees on
  $k$ vertices having degree $d$. Setting $k=1$, we get the degree
  sequence.
\end{proof}

\section{Edge-subgraph posets and edge reconstruction}
\label{sec-esp}

In this section, we classify graphs that are not $Q$-reconstructible.
We show that if the edge reconstruction conjecture is true, then graphs
that are not $Q$-reconstructible, except finitely many, have a simple
structure: if $G$ and $H$ are distinct unlabelled graphs that have
isomorphic edge subgraph poset, then except in finitely many cases,
$\{G,H\} = \{K_{1,m},mK_2\}$ for some $m \geq 2$, or
$\{G,H\} = \{pK_3+qK_{1,3}+F, qK_3+pK_{1,3}+F\}$, where $p\neq q$ and
$F$ itself is a graph with quite simple structure.

\begin{figure}[ht]
  \begin{center}
    \includegraphics[scale=1.0]{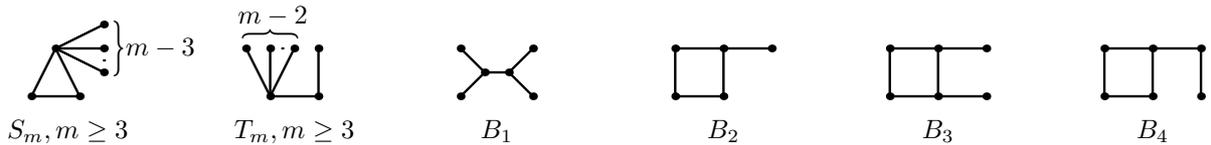}
  \end{center}
  \vspace*{-2ex}
  \caption[]{Some graphs referenced in Theorem~\ref{thm-esp}}
  \label{fig-names}
\end{figure}

\begin{notation} Graphs $B_1,\ldots, B_4$, and $m$-edge graphs $S_m$ and
  $T_m$, where $(m\geq 3)$, which are frequently referred to in the
  proofs, are as shown in Figure~\ref{fig-names}. %
  \nomenclature[4b]{$B_1,\ldots, B_4$}{certain special graphs} %
  \nomenclature[4s]{$S_m,m \geq 3$}{certain special graphs} %
  \nomenclature[4t]{$T_m,m \geq 3$}{certain special graphs} %
  We denote by $\eminus{K_4}$ the isomorphism class of a graph obtained
  by deleting an edge from a copy of $K_4$. %
  Let
  \[
  \mcf\coloneqq\{P_n\mid n\geq 2 \}\bigcup\{C_n,n\geq 4\}\bigcup \{S_4,
  \eminus{K_4}, K_4\},\text{ and }
  \]
  Let $\nn^{(\mcf)}$ be the set of all finite unlabelled graphs
  (including the null graph) with components from $\mcf$. %
  \nomenclature[4nf]{$\nn^{(\mcf)}$}{set of all unlabelled graphs with
    components from $\mcf$} %

  For an unlabelled graph $F$, define $\eplus{F}$ to be the set of
  unlabelled graphs that can be obtained by adding a new edge to a copy
  of $F$, where the added edge may have 0, 1, or 2 end-vertices in the
  copy of $F$; formally,
  $\eplus{F}\coloneqq\{\iso{H_1\cup H_2}\mid H_1\in F,H_2\in K_2\}$.
  \nomenclature[4fe]{$\eplus{F}$}{set of unlabelled graphs obtained by
    adding an edge to $F$} %
  For example, $\eplus{K_{1,3}} = \{K_{1,3}+K_2, K_{1,4},T_4,S_4\}$.
\end{notation}

\begin{thm} \label{thm-esp}\leavevmode %
  \makeatletter\renewcommand\p@enumi{\thethm--}\makeatother
  \makeatletter\renewcommand\p@enumii{\thethm--\theenumi}\makeatother
  \begin{enumerate}%
  \item\label{A} The graphs in each of the following sets have the same
    abstract edge-subgraph poset:
    \begin{enumerate}%
    \item \label{6e} %
      $\{K_3, K_{1,3}, 3K_2\}$, %
      $\{P_4, K_{1,2}+K_2\}$, %
      $\{P_4+K_2, T_4\}$, %
      $\{C_4, 2K_{1,2}\}$, %
      $\{C_4 + K_2, B_1\}$, %
      $\{P_6, B_2\}$, and %
      $\{B_3, B_4\}$;
    \item \label{2e} $\{K_{1,m}, mK_2\}$, for all $m > 1$;
    \item \label{7e} $\{pK_3+qK_{1,3}+F, qK_3+pK_{1,3}+F\}$, where
      $p \neq q$ and $F \in \nn^{(\mcf)}$.%
    \end{enumerate}
  \item\label{B}The edge reconstruction conjecture is true if and only
    if all graphs, except the graphs listed above, are
    $Q$-reconstructible.
  \end{enumerate}
\end{thm}

\begin{lem} \label{lem-1c-closed} If $G \in pK_3+qK_{1,3}+F$ for some
  $p,q \in \nn$ and $F \in \nn^{(\mcf)}$, and $G_i$ is an edge-subgraph
  of $G$, then $G_i \in p_iK_3+q_iK_{1,3}+F_i$ for some
  $p_i,q_i \in \nn$ and $F_i \in \nn^{(\mcf)}$.
\end{lem}
\begin{proof}
  The claim is proved by verifying that all proper edge-subgraphs of
  $K_3$, all proper edge-subgraphs of $K_{1,3}$, and all edge-subgraphs
  of each graph in $\mcf$ are in $\nn^{(\mcf)}$.
\end{proof}

In the following, we say that labelled graphs $G$ and $H$ are {\em
  conjugates} if there exist $p,q \in \nn$ (possibly equal) and
$F \in \nn^{(\mcf)}$ such that $G \in pK_3 + qK_{1,3} + F$ and
$H \in pK_{1,3} + qK_3 + F$.

\begin{lem}\label{lem-1c-weights1}
  Let the graphs $G$ and $H$ be conjugates. Then there exists a
  bijection $f\colon 2^{E(G)}\to 2^{E(H)}$ such that for all
  $E\subseteq E(G)$, the edge-subgraphs $G[E]$ and $H[f(E)]$ are
  conjugates.
\end{lem}

\begin{proof} %
  Order the components $G^i$ of $G$ and components $H^i$ of $H$ such
  that $G^i \in K_3$ if and only if $H^i\in K_{1,3}$, and
  $G^i \in K_{1,3}$ if and only if $H^i\in K_3$, and $G^i\cong H^i$
  otherwise.

  For each $i \in \{1,\ldots, p+q+r\}$, define a bijection
  $f_i \colon 2^{E(G^i)}\bij 2^{E(H^i)}$ such that
  \begin{equation}
    \label{eq-map-1c} 
    H^i[f_i(E)]
    \begin{cases}
      \in K_{1,3} & \text{ if } G^i[E]\in K_3, \\
      \in K_3 &  \text{ if } G^i[E]\in K_{1,3}, \\
      \cong G^i[E] & \text{ otherwise}.
    \end{cases}
  \end{equation}
  Such a bijection always exists due to the chosen ordering of the
  components of $G$ and $H$, and since for each graph $F \in \mcf$, we
  have $\esub{K_3}{F} = \esub{K_{1,3}}{F}$; and if $G^i\in K_3$ or
  $G^i\in K_{1,3}$, then any bijection $f_i$ such that $|f_i(E)| = |E|$
  serves the purpose.

  Define $f\colon 2^{E(G)}\to 2^{E(H)}$ by extending the component-wise
  maps $f_i$ such that for all $E\subseteq E(G)$,
  \begin{equation*}
      f(E) = \bigcup_i f_i(E\cap E(G^i)).
  \end{equation*}
  The bijection $f$ has the desired property.
\end{proof}

\begin{lem}\label{lem-1c-weights2}
  Let $G^i, H^i$ and $G^j, H^j$ be pairs of conjugates. Then
  $\esub{G^i}{G^j} = \esub{H^i}{H^j}$.
\end{lem}

\begin{proof}
  Let $f\colon 2^{E(G^j)} \to 2^{E(H^j)}$ be a bijection as defined in
  the statement of Lemma~\ref{lem-1c-weights1}. Let $E\subseteq
  E(G^j)$.
  We have $G^j[E] \cong G^i$ if and only if $H^j[f(E)] \cong H^i$; hence
  $\esub{G^i}{G^j} = \esub{H^i}{H^j}$.
\end{proof}

\begin{proof}[Proof of Theorem~\ref{A}]
  For the graphs listed in Theorem~\ref{6e} and~\ref{2e}, the claim is
  proved by constructing the abstract edge-subgraph poset for each graph
  and verifying that the graphs within each set have the same abstract
  edge-subgraph poset.

  Let $G = pK_3 + qK_{1,3} + F$ and $H = pK_{1,3} + qK_3 + F$, where
  $p,q \in \nn$, $p \neq q$, and $F \in \nn^{(\mcf)}$.
  Lemmas~\ref{lem-1c-closed},~\ref{lem-1c-weights1},
  and~\ref{lem-1c-weights2} imply that there exists an isomorphism from
  $\esp{G}$ to $\esp{H}$ that maps an edge-subgraph
  $p_iK_3+q_iK_{1,3}+F_i$ of $G$ to an edge-subgraph
  $p_iK_{1,3}+q_iK_3+F_i$ of $H$, where $p_i,q_i \in \nn$ and
  $F_i \in \nn^{(\mcf)}$. But $G \neq H$ since $p\neq q$. This proves the
  result for graphs listed in Theorem~\ref{7e}.
\end{proof}

\begin{proof}[Proof of Theorem~\ref{B} (the `if' part)] The assertion
  follows from the following facts.
\begin{enumerate}%
\item All the graphs listed in Theorem~\ref{A} that have at least 4
  edges are edge reconstructible since disconnected graphs on 4 or more
  edges, trees on 4 or more edges, and unicyclic graphs on 4 or more
  edges are edge reconstructible; see \citet{kelly1957},
  \citet{greenwell.hemminger-1969}, and \citet{manvel-1969a}.

\item All $Q$-reconstructible graphs are edge reconstructible since the
  abstract edge-subgraph poset of a graph can be constructed given its
  edge-deck.  \qedhere
\end{enumerate}
\end{proof}
\smallskip
\begin{proof}[Proof of Theorem~\ref{B} (outline of the `only if' part)]
  The `only if' part is proved in 3 steps.
  \begin{enumerate}%
  \item In Proposition~\ref{p0-th}, the result is proved for graphs with
    at least 4 edges, containing both $K_{1,2}+ 2K_2$ and $T_4$ as
    subgraphs.
  \item In Proposition~\ref{p1-th}, a proof of the result is sketched
    for graphs with at most 7 edges.
  \item In Proposition~\ref{p2-th}, the result is proved for graphs with
    at least 7 edges that do not have at least one of the graphs
    $K_{1,2}+ 2K_2$ and $T_4$ as a subgraph. \qedhere %
  \end{enumerate}
\end{proof}
In the rest of this section, we prove
Propositions~\ref{p0-th},~\ref{p1-th} and~\ref{p2-th}, followed by some
corollaries and questions.

A {\em legitimate labelling} of $\esp{G}$ is a one-to-one map
$\pi\colon \esp{G} \to \unlab$ such that $\pi(g_1) = K_2$ and
$\we(g_i,g_j) = \esub{\pi(g_i)}{\pi(g_j)}$ for all
$g_i,g_j \in \esp{G}$. We say that an element $g_i$ (or a subset $S$) of
$\esp{G}$ is {\em uniquely labelled} if for all legitimate labelling
maps $\pi$ and $\sigma$ from $\esp{G}$ to $\unlab$, we have
$\pi(g_i) = \sigma(g_i)$ (or $\pi(g_i) = \sigma(g_i)$ for all $g_i$ in
$S$). Thus $G$ is $Q$-reconstructible if and only if the maximal element
of $\esp{G}$ is uniquely labelled. The $Q$-reconstructibility of $G$
does not imply that $\esp{G}$ is uniquely labelled; for example, $S_4$
is $Q$-reconstructible, but given a legitimate labelling of $\esp{S_4}$,
the labels $K_{1,3}$ and $K_3$ may be interchanged keeping all other
labels fixed to obtain another legitimate labelling.
\begin{lem}
  \label{lem1.1} The graph $K_{1,2}+ 2K_2$ is $Q$-reconstructible, and
  the abstract edge-subgraph poset of $K_{1,2}+ 2K_2$ is uniquely
  labelled.
\end{lem}

\begin{proof} 
  In Figure~\ref{fig-part-labelled}, we have
  $\esp{g_9} \cong \espu{K_{1,2}+ 2K_2}$. Let $\pi$ be a legitimate
  labelling of $\esp{g_9}$. We show that $\pi$ is unique.  If
  $\pi(g_2) = K_{1,2}$ and $\pi(g_3) = 2K_2$, then $\pi(g_5) = P_4$, and
  $\pi(g_4) \in \{K_{1,3},K_3\}$. Such a labelling cannot be extended so
  as to assign a legitimate label to $g_9$, since it would imply that
  $\pi(g_9)$ has edge-deck $\{P_4,P_4,K_3,K_3\}$ or
  $\{P_4,P_4,K_{1,3},K_{1,3}\}$. But neither $\{P_4,P_4,K_3,K_3\}$ nor
  $\{P_4,P_4,K_{1,3},K_{1,3}\}$ is a legitimate edge-deck. Hence $g_9$
  does not have a legitimate labelling $\pi$ in which
  $\pi(g_2) = K_{1,2}$ and $\pi(g_3) = 2K_2$. %
  If $\pi(g_2) = 2K_2$ and $\pi(g_3) = K_{1,2}$, then the edge
  reconstructibility of each of $3K_2$, $K_{1,2}+K_2$, and
  $K_{1,2}+2K_2$ determines that $\pi(g_4) = 3K_2$,
  $\pi(g_5) = K_{1,2}+K_2$, and $\pi(g_9) = K_{1,2}+2K_2$.
\end{proof}

\begin{prop}
  \label{p0-th}
  If the edge reconstruction conjecture is true, then all graphs with at
  least 4 edges that contain both $K_{1,2}+ 2K_2$ and $T_4$ as subgraphs
  are $Q$-reconstructible. %
\end{prop}

\begin{proof} Given $\esp{G}$, where $\epsilon(G) \geq 4$, we recognise
  that $K_{1,2}+ 2K_2$ and $T_4$ are subgraphs of $G$; then we uniquely
  label all elements of $\esp{G}$ of rank 3; then, assuming that the
  edge reconstruction conjecture is true, we extend the unique labelling
  to all of $\esp{G}$, thereby implying the $Q$-reconstructibility of
  $G$.

  {\em Recognising that $K_{1,2}+ 2K_2$ is a subgraph of $G$:} By
  Lemma~\ref{lem1.1}, $K_{1,2}+ 2K_2$ is a subgraph of $G$ if and only
  if there exists $g_i\in \esp{G}$ such that
  $\esp{g_i}\cong \espu{K_{1,2}+ 2K_2}$; and if such a $g_i$ exists, then
  $g_i$ is unique and $\esp{g_i}$ is uniquely labelled. Suppose that
  there exists such a $g_i$. Let $\pi$ be an arbitrary labelling of
  $\esp{G}$; it is unique on $\esp{g_i}$. Let
  $\pi(g_1) = K_2$, $\pi(g_2) = 2K_2$, $\pi(g_3) = K_{1,2}$,
  $\pi(g_4) = 3K_2$ and $\pi(g_5) = K_{1,2}+K_2$.

  {\em Recognizing that $T_4$ is a subgraph of $G$ given that
    $K_{1,2}+ 2K_2$ is a subgraph of $G$:} We have
  $\espu{T_4} \cong \espu{P_4+K_2}$, and there is no other graph $H$
  such that $\espu{H} \cong \espu{T_4}$. Let $g_j\in \esp{G}$ such that
  $\esp{g_j} \cong \espu{T_4}$. Now $3K_2$ is a subgraph of $P_4+K_2$
  but not of $T_4$, hence $\pi(g_j)=P_4+K_2$ if $g_4 \epo g_j$, else
  $\pi(g_j)=T_4$. Suppose that it is the latter case.

  {\em Uniquely labelling all elements of rank $3$:} The unique partial
  labelling constructed above extends to all elements of $\esp{G}$ of
  rank $3$ as follows. We have
  $\espu{3K_2} \cong \espu{K_{1,3}} \cong \espu{K_3}$, but the labels
  $3K_2$ and $T_4$ have been assigned uniquely, and $K_{1,3}$ is a
  subgraph of $T_4$ while $K_3$ is not a subgraph of $T_4$; thus the
  labels $K_{1,3}$ and $K_3$ are uniquely assigned. Other graphs with 3
  edges are uniquely labelled since they are edge reconstructible, and
  graphs with 2 edges are uniquely labelled.

  {\em Reconstruction:} If the edge reconstruction conjecture is true,
  then for all $k\geq 3$, a unique labelling of all elements of rank $k$
  extends to a unique labelling of all elements of rank $k+1$. Since all
  elements of rank 3 are uniquely labelled, an induction on the number
  of edges implies the result.
\end{proof}

\begin{lem}
\label{lem1.3}
For all $m \geq 4$, if $G \in \eplus{K_{1,m}}\setminus \{K_{1,m+1}\}$,
then $G$ is $Q$-reconstructible, and the label $K_{1,3}$ and the label
$K_3$ (if $K_3$ is a subgraph of $G$) are uniquely assigned to elements
of $\esp{G}$.
\end{lem}

\begin{proof} We have
  $\eplus{K_{1,m}} = \{K_{1,m+1}, S_{m+1}, T_{m+1}, K_{1,m}+K_2\}$,
  where the graphs $S_{m+1}$ and $T_{m+1}$ are as shown in
  Figure~\ref{fig-names}. The affirmation is true for $m = 4$, i.e., the
  graphs $S_5$, $T_5$ and $K_{1,4}+K_2$ are $Q$-reconstructible. Hence
  we assume that $m > 4$. Now $G = S_{m+1} $ (or $G = T_{m+1}$, or
  $G = K_{1,m}+K_2$) if and only if the following conditions are
  satisfied:
  \begin{enumerate}
  \item $\esp{G}$ is not totally ordered
  \item there exists $g_i \in \esp{G}$ such that $\rho(g_i) = m$, and
    $\esp{g_i}$ is totally ordered
  \item there exists $g_j \in \esp{G}$ such that
    $\esp{g_j} \cong \espu{S_{5}}$ (or %
    $\esp{g_j} \cong \espu{T_{5}}$ or %
    $\esp{g_j} \cong \espu{K_{1,4}+K_2}$, respectively).
  \end{enumerate}
  The necessity is directly verified for each graph in
  $\eplus{K_{1,m}}\setminus \{K_{1,m+1}\}$. %
  For sufficiency, the first two conditions imply that either
  $G \in\eplus{K_{1,m}}\setminus \{K_{1,m+1}\}$ or
  $G \in \eplus{(mK_2)}\setminus \{(m+1)K_2\}$; %
  the third condition implies that
  $G \not \in \eplus{(mK_2)}\setminus \{(m+1)K_2\}$; %
  finally, $S_5$ is a subgraph of $S_{m+1}$ but not of $T_{m+1}$ or
  $K_{1,m}+K_2$ (and $T_5$ is a subgraph of $T_{m+1}$ but not of
  $S_{m+1}$ or $K_{1,m}+K_2$, and $K_{1,4}+K_2$ is a subgraph of
  $K_{1,m}+K_2$ but not of $T_{m+1}$ or $S_{m+1}$).

  \smallskip
  
  {\em Assigning the labels $K_{1,3}$ and $K_3$:} If $G = S_{m+1}$, then
  $K_3$ and $K_{1,3}$ are subgraphs of $G$, and $3K_2$ is not a subgraph
  of $G$. But $\esub{K_3}{S_{m+1}} = 1$ and
  $\esub{K_{1,3}}{S_{m+1}} > 1$. Hence the labels $K_{1,3}$ and $K_3$
  are uniquely assigned.
  If $G = T_{m+1}$ or $G = K_{1,m}+K_2$, then neither $K_3$ nor $3K_2$
  is a subgraph of $G$, hence the label $K_{1,3}$ is uniquely assigned.
\end{proof}

\begin{lem}
  \label{lem1.4}
  For all $m \geq 4$, if
  $G \in \eplus{(mK_2)}\setminus \{(m+1)K_2\}$, then $G$ is
  $Q$-reconstructible. A unique element in $G$ is labelled $3K_3$.
\end{lem}
\begin{proof} We have
  $\eplus{(mK_2)} = \{(m+1)K_2, K_{1,2} + (m-1)K_2, P_4+(m-2)K_2 \}$ for
  $m \geq 2$.  The affirmation is true for $m = 4$, i.e., the graphs
  $K_{1,2}+3K_2$ and $P_4+2K_2$ are $Q$-reconstructible. Hence we assume
  that $m > 4$.  Now $G = K_{1,2} + (m-1)K_2$ (or $G = P_4+(m-2)K_2$) if
  and only if the following conditions are satisfied:
  \begin{enumerate}
  \item $\esp{G}$ is not totally ordered
  \item there exists $g_i \in \esp{G}$ such that $\rho(g_i) = m$, and
    $\esp{g_i}$ is totally ordered
  \item there exists $g_j \in \esp{G}$ such that
    $\esp{g_j} \cong \espu{K_{1,2}+3K_2}$ (or, respectively,
    $\esp{g_j} \cong \espu{P_4+2K_2}$).
  \end{enumerate}
  The necessity is directly verified for each of the graphs
  $K_{1,2} + (m-1)K_2$ and $P_4+(m-2)K_2 $. %
  For sufficiency, the first two conditions imply that either
  $G \in\eplus{K_{1,m}}\setminus \{K_{1,m+1}\}$ or
  $G \in \eplus{(mK_2)}\setminus \{(m+1)K_2\}$; %
  the third condition implies that
  $G \not \in \eplus{K_{1,m}}\setminus \{K_{1,m+1}\}$, hence $g_i$ must
  be labelled $mK_2$; %
  finally, $G = K_{1,2}+(m-1)K_2$ if
  $\esub{mK_2}{G} = \we(g_i,g_{\scriptscriptstyle M}) = 2$, and
  $G = P_4+(m-1)K_2$ if
  $\esub{mK_2}{G} = \we(g_i,g_{\scriptscriptstyle M}) = 1$.

  \smallskip

  {\em Assigning the label $3K_2$:} Neither $K_3$ nor $K_{1,3}$ is a
  subgraph of any graph $G$ in $\eplus{(mK_2)}\setminus \{(m+1)K_2\}$,
  hence the label $3K_2$ is uniquely assigned.
\end{proof}

\begin{prop}
  \label{p1-th}
  Graphs with at most 7 edges, except the ones listed in
  Theorem~\ref{A}, are $Q$-reconstructible.
\end{prop}

\begin{proof}[A sketch of the proof]
  The proof requires looking at several straightforward
  cases. Therefore, we only indicate two techniques, besides Lemmas and
  Propositions~\ref{lem1.1} to~\ref{lem1.4}, that we use to prove the
  result efficiently.

  Let $\{G,H\}$ be a $Q$-pair of graphs with $m$ edges, where
  $m \geq 4$. Let $f$ be an isomorphism from $\espu{G}$ to
  $\espu{H}$. Then there exists a $Q$-pair $\{G_i,H_i\}$ of edge-deleted
  subgraphs $G_i\epo G$ and $H_i\epo H$, such that $f(G_i) = H_i$;
  otherwise the edge reconstruction conjecture would imply that $G$ and
  $H$ are isomorphic. Thus once all $Q$-pairs $\{G_i, H_i\}$ of graphs
  with $m-1$ edges are enumerated, we only need to consider pairs
  $\{G,H\}$ such that $G\in\eplus{G_i}$ and $H\in\eplus{H_i}$ as
  probable candidates for $Q$-pairs on $m$ edges. Many graphs in the
  sets $\eplus{G_i}$ and $\eplus{H_i}$ are proved to be
  $Q$-reconstructible using Lemmas~\ref{lem1.1} to~\ref{lem1.4}. This
  significantly reduces the number of cases that need to be analysed.

  The second technique is assigning unique labels  to some elements of
  $\esp{G}$. If $G_i$ is a $Q$-reconstructible graph, and
  $\esp{g_i}\cong\espu{G_i}$, then only $g_i$ can be assigned the label
  $G_i$. Such an assignment of a label may uniquify labels $G_j$ of
  graphs which may not otherwise be $Q$-reconstructible. We used this
  idea in Lemmas~\ref{lem1.1} to~\ref{lem1.4}; for example, presence of
  a subgraph $K_{1,2}+ 2K_2$, which is $Q$-reconstructible, uniquifies
  labels $K_{1,2}, 2K_2, 3K_2, K_{1,2}+K_2$, even though neither of
  these graphs is $Q$-reconstructible; similarly, if a graph in
  $\eplus{K_{1,4}}$ is a subgraph of $G$, then it uniquifies labels
  $K_{1,3}, K_3, 3K_2$ even though neither of these graphs is
  $Q$-reconstructible. Once many labels have been fixed, the
  $Q$-reconstructibility may follow quickly.
\end{proof}

Let the graphs listed in Theorem~\ref{A} be grouped into 3 families
defined below:
\begin{align*}
  \mcc_1 & \coloneqq \{P_4, K_{1,2}+K_2, P_4+K_2, T_4, C_4, 2K_{1,2},
           C_4 + K_2, B_1, P_6, B_2, B_3, B_4\},\\
  \mcc_2 & \coloneqq \{K_{1,m}\mid m > 1 \} \bigcup \{mK_2 \mid m > 1\},
           \text{ and }\\
  \mcc_3 & \coloneqq \{pK_3+qK_{1,3}+F\mid p\neq q \text{ and } F \in
           \nn^{(\mcf)}\}.
\end{align*}

\begin{prop}
\label{p2-th}
If the edge reconstruction conjecture is true, then graphs with at least
7 edges, except the ones listed in Theorem~\ref{A}, are
$Q$-reconstructible.
\end{prop}

\begin{proof}
  We assume that the edge reconstruction conjecture is true.  Let $G$ be
  a graph to be $Q$-reconstructed and $\epsilon(G) \geq 7$. %
  By Propositions~\ref{p0-th}, we assume that $G$ does not contain at
  least one of the graphs $K_{1,2}+ 2K_2$ and $T_4$ as a subgraph. %
  The graphs in $\mcc_1$ have 6 or fewer edges. Hence we show by
  induction on $\epsilon(G)$, that if $G$ is not $Q$-reconstructible,
  then $G\in \mcc_2\cup\mcc_3$. %

  We take $\epsilon(G) = 7$ as the base case, for which the result
  follows from Theorem~\ref{A} and Proposition~\ref{p1-th} since the
  graphs in $\mcc_1$ all have at most 6 edges. Suppose that the
  affirmation is true when $7 \leq \epsilon(G) \leq m$. Let $G$ be an
  $(m+1)$-edge graph that is not $Q$-reconstructible. We show that
  $G \in \mcc_2\cup\mcc_3$.

  There must exist an $m$-edge subgraph $G_i$ of $G$ that is not
  $Q$-reconstructible; otherwise we would be able to construct the
  edge-deck of $G$ from its abstract edge-subgraph poset, and then the
  edge reconstruction conjecture would imply the $Q$-reconstructibility
  of $G$. By induction hypothesis, each $m$-edge subgraph of $G$ that is
  not $Q$-reconstructible is in $\mcc_2\cup \mcc_3$.

  \begin{pclaims}

    \pclaim\label{claim-gj-in-8b} Let $G_i$ be an $m$-edge subgraph $G$
    that is not $Q$-reconstructible. If $G_i \in \mcc_2$, then
    $G \in \mcc_2$.
    
    \ppf If $G_i = K_{1,m}$, then $G\in\eplus{K_{1,m}}$. By
    Lemma~\ref{lem1.3}, $G = K_{1,m+1}\in\mcc_2$. If $G_i = mK_2$, then
    $G\in\eplus{(mK_2)}$. By Lemma~\ref{lem1.4},
    $G=(m+1)K_2\in\mcc_2$. \pqed

    Therefore, in the following, we assume that all $m$-edge subgraphs
    of $G$ that are not $Q$-reconstructible are in $\mcc_3$, and show
    that $G$ is in $\mcc_3$.

    \pclaim \label{claim-no-t1} The graph $K_{1,2}+ 2K_2$ is a
    subgraph of $G$, and $T_4$ is not a subgraph of $G$.

    \ppf Let $G_i$ be an $m$-edge subgraph of $G$ that is not
    $Q$-reconstructible. The assumptions $G_i \in \mcc_3$ and
    $e(G_i) \geq 7$ imply that $K_{1,2}+ 2K_2 \epo G_i \epo G$. But it
    follows from Proposition~\ref{p0-th} that both $K_{1,2}+ 2K_2$ and
    $T_4$ cannot be subgraphs of $G$, which implies the claim. \pqed

    \pclaim \label{claim-components} The parameters
    $\esub{S_4}{G}, \esub{\eminus{K_4}}{G}, \text{ and } \esub{K_4}{G}$
    are $Q$-reconstructible.
    The graphs $S_4$, $\eminus{K_4}$ and $K_4$ can appear as
    subgraphs of $G$ only within components of $G$ on 4 vertices.

    \ppf The graphs $S_4$, $\eminus{K_4}$ and $K_4$ are
    $Q$-reconstructible, implying the first part. If any of the graphs
    $S_4$, $\eminus{K_4}$ and $K_4$ is in a component of $G$ on 5 or
    more vertices, then $T_4$ is a subgraph of $G$; but $T_4$ has been
    eliminated by Claim~\ref{claim-no-t1}. \pqed

    \pclaim \label{claim-counting-T} All elements of $\esp{G}$ of rank 3
    and 4, except possibly the ones corresponding to $K_{1,3}$, $K_3$,
    $K_{1,3}+K_2$ and $K_3+K_2$, are uniquely labelled. The parameters
    $\esub{K_3}{G}+\esub{K_{1,3}}{G}$, and
    $\esub{K_3+K_2}{G}+\esub{K_{1,3}+K_2}{G}$, and $\esub{T}{G}$, are
    $Q$-reconstructible, where $T$ is any graph with at most 4 edges.

    \ppf The claim follows from the following two statements. By
    Claim~\ref{claim-no-t1}, $K_{1,2}+ 2K_2$ is a subgraph of $G$. By
    Lemma~\ref{lem1.1}, $K_{1,2}+ 2K_2$ is $Q$-reconstructible, and the
    elements of $\espu{K_{1,2}+ 2K_2}$ are uniquely labelled (in
    particular, a unique element is assigned the label $3K_2$).  \pqed

    \pclaim \label{claim-K4} If $G$ contains $K_4$, then $G \in \mcc_3$.

    \ppf Since $K_4$ is $Q$-reconstructible, $G$ contains $K_4$ if and
    only if there is a unique element $g_j\in\esp{G}$ such that
    $\esp{g_j} \cong \espu{K_4}$. Assume that that is the case. By
    Claim~\ref{claim-components}, $K_4$ can only occur as a component of
    $G$.
    Let $g_k$ be an element of $\esp{G}$ such that $\rho(g_k) = m$, and
    $\we(g_j,g_k) \neq \we(g_j,g_M)$; i.e.,
    $\esub{K_4}{G_k} \neq \esub{K_4}{G}$. Such an element must exist
    since deleting an edge from a component isomorphic to $K_4$ gives a
    graph such as $G_k$. The graph $G_k$ cannot be $Q$-reconstructible;
    otherwise $G$ would be obtained by adding an edge to any component
    of $G_k$ that is isomorphic to $\eminus{K_4}$, and hence $G$ would
    be $Q$-reconstructible as well. Hence $G_k\in \mcc_3$; and adding an
    edge to a component of $G_k$ that is isomorphic to $\eminus{K_4}$
    results in a graph in $\mcc_3$.  \pqed

    Therefore, we assume that $K_4$ is not a subgraph of $G$.

    \pclaim \label{claim-K4-e} If $G$ contains $\eminus{K_4}$, then
    $G\in\mcc_3$.

    \ppf Since $\eminus{K_4}$ is $Q$-reconstructible, $G$ contains
    $\eminus{K_4}$ if and only if there is a unique element
    $g_j\in\esp{G}$ such that $\esp{g_j} \cong
    \espu{\eminus{K_4}}$.
    Assume that that is the case. Since we have assumed that $G$ does
    not contain $K_4$, by Claim ~\ref{claim-components}, the graph
    $\eminus{K_4}$ can only occur as a component of $G$.
    Let $g_k$ be an element of $\esp{G}$ such that $\rho(g_k) = m$, and
    $\esub{\eminus{K_4}}{G_k} = \esub{\eminus{K_4}}{G}-1$, and
    $\esub{K_3}{G_k}+\esub{K_{1,3}}{G_k} =
    \esub{K_3}{G}+\esub{K_{1,3}}{G}-4$.
    Such an element $g_k$ must exist since there is a unique edge in
    $\eminus{K_4}$ that belongs to two 3-stars and two triangles.  Now,
    as in Claim~\ref{claim-K4}, the graph $G_k$ cannot be
    $Q$-reconstructible; otherwise $G$ would be obtained by adding an
    edge to any component of $G_k$ that is isomorphic to a $4$-cycle,
    and hence $G$ would be $Q$-reconstructible as well. Hence
    $G_k\in \mcc_3$; and adding an edge to a component $4$-cycle of
    $G_k$ results in a graph in $\mcc_3$. \pqed
  
    Therefore, we assume that $\eminus{K_4}$ is not a subgraph of $G$.

    \pclaim \label{claim-B1} If $S_4$ is a component of $G$, then
    $G\in\mcc_3$.

    \ppf Since $S_4$ is $Q$-reconstructible, we recognise that $S_4$ is
    a component of $G$ if and only if there exists $g_j\in\esp{G}$ such
    that $\esp{g_j} \cong \espu{S_4}$. Assume that that is the
    case. Since we have assumed that $G$ does not contain $K_4$ or
    $\eminus{K_4}$, the graph $S_4$ can only occur as a component of
    $G$.
    Claim~\ref{claim-counting-T} implies that there is a unique
    $g_k \in \esp{G}$ that is labelled $P_4$ (even though $P_4$ itself
    is not $Q$-reconstructible).
    There is an edge in $S_4$ that belongs to exactly one $P_4$,
    therefore, there exists $g_\ell \in \esp{G}$ of rank $m$ such that
    $\we(g_j,g_\ell) = \we(g_j,g_M) -1$, and
    $\we(g_k,g_\ell) = \we(g_k,g_M) -1$.  The two conditions mean,
    respectively, $\esub{S_4}{G_\ell} = \esub{S_4}{G}-1$ and
    $\esub{P_4}{G_\ell} = \esub{P_4}{G}-1$.
    Hence $G$ is obtained from $G_\ell$ by adding an edge in a component
    isomorphic to $P_4$ so as to create a component isomorphic to
    $S_4$. As in Claims~\ref{claim-K4} and~\ref{claim-K4-e}, the graph
    $G_\ell$ cannot be $Q$-reconstructible; otherwise $G$ would be
    $Q$-reconstructible as well. Hence then $G_\ell \in \mcc_3$,
    implying that $G \in \mcc_3$.  \pqed
  
    Therefore, we assume that $S_4$ is not a subgraph of $G$.

    \pclaim If $G$ contains a cycle, then the cycle is a component. If
    $G$ contains a path on 4 or more vertices, then the path is either a
    component or a subgraph of a component that is either a cycle or a
    path.

    \ppf Both parts follows from the assumption that neither $S_4$ nor
    $T_4$ is a subgraph of $G$. \pqed

    \pclaim \label{claim-pn-cn} All elements of $\esp{G}$ corresponding
    to paths and cycles, except possibly $K_3$, are uniquely labelled.

    \ppf By Claim~\ref{claim-counting-T}, the elements of $\esp{G}$
    corresponding to $P_2$, $P_3$, and $P_4$ are uniquely labelled. For
    $n \geq 4$, the graphs $P_n$ and $C_n$ are edge reconstructible. Now
    the claim is proved by induction on $n\geq 4$. \pqed

    Therefore, we assume that $G$ itself is not a path or a cycle.

    \pclaim \label{claim-pn} If $P_n, n \geq 2$ is a component of
    $G$, then $G\in\mcc_3$.

    \ppf Let $g_j$ be an element of rank $m$ in $\esp{G}$. The graph
    $G_j$ is obtained from $G$ by deleting an edge at the end of a
    component path $P_n$ if and only if
    $\esub{P_k}{G}=\esub{P_k}{G_j}+1$ for all $k\leq n$, and
    $\esub{P_{n+1}}{G}=\esub{P_{n+1}}{G_j}$. These conditions are
    recognised from $\esp{G}$ by Claim~\ref{claim-pn-cn}.  Assume that
    $g_j$ is such an element; hence $G$ is obtained by adding an edge at
    the end of a path isomorphic to $P_{n-1}$ in $G_j$.
    The graph $G_j$ is not $Q$-reconstructible, since otherwise $G$
    would be $Q$-reconstructible also. Hence both $G_j$ and $G$ must be
    in $\mcc_3$. \pqed

    Therefore, we assume that $G$ does not contain a path on 2 or more
    vertices as a component.

    \pclaim If $G$ has a component isomorphic to a cycle
    $C_n, n \geq 4$, then $G\in\mcc_3$.

    \ppf Let $g_j$ be an element of rank $m$ in $\esp{G}$. The graph
    $G_j$ is obtained from $G$ by deleting an edge of a component cycle
    $C_n$ if and only if $\esub{C_n}{G_j}= \esub{C_n}{G}-1$. This
    condition can be recognised from $\esp{G}$ by
    Claim~\ref{claim-pn-cn}. Assume that $g_j$ is such an element; hence
    $G$ is obtained by adding an edge to $G_j$ joining the end vertices
    of a component $P_n$ in $G_j$.
    As in earlier claims, the graph $G_j$ cannot be $Q$-reconstructible;
    otherwise $G$ would be $Q$-reconstructible also. Hence both $G_j$
    and $G$ must be in $\mcc_3$.
  \end{pclaims}
  
  Therefore, we assume that $G$ does not contain a cycle on 4 or more
  vertices. Now the only remaining graphs are the graphs in which all
  components are isomorphic to $K_3$ or $K_{1,3}$, and their total
  number is $Q$-reconstructible by Claim~\ref{claim-counting-T},
  completing the proof of Proposition~\ref{p2-th}.
\end{proof}

The method of the proof of Proposition~\ref{p2-th} may be applied to any
class of graphs that is closed under edge-deletion.

\begin{cor}
\label{cor-trees}
Acyclic graphs (i.e., trees and forests) that are not in class
$\mcc_1\cup\mcc_2\cup\mcc_3$ are $Q$-reconstructible.
\end{cor}

\begin{proof} Acyclic graphs with four or more edges are edge
  reconstructible. The class of acyclic graphs is closed under
  edge-deletion. Hence the proof of Proposition~\ref{p2-th} may be
  restricted to the class of acyclic graphs to show that acyclic graphs
  that are not $Q$-reconstructible belong to
  $\mcc_1\cup\mcc_2\cup\mcc_3$.
\end{proof}

\begin{prop}[\citet{muller1977}] \label{prop-muller} All graphs $G$ such
  that $2^{\epsilon(G)-1} > \nu(G)!$ are edge reconstructible.
\end{prop}

The following result is a weaker version of M\"uller's result for the
$Q$-reconstruction problem.

\begin{cor}\label{cor-muller}
  The edge reconstruction conjecture is true if and only if all graphs
  $G$ such that $2^{\epsilon(G)-1}\leq \nu(G)!$, except the ones in the
  class $\mcc_1\cup\mcc_2\cup\mcc_3$, are $Q$-reconstructible.
\end{cor}

\begin{proof} Theorem~\ref{thm-esp} and Proposition~\ref{prop-muller}
  imply the result.
\end{proof}

\begin{rem}
  Corollary~\ref{cor-muller} implies that if all graphs $G$ such that
  $G\not\in \mcc_1\cup\mcc_2\cup\mcc_3$ and
  $2^{\epsilon(G)-1}\leq \nu(G)!$ are $Q$-reconstructible, then graphs
  $G$ such that $G\not\in \mcc_1\cup\mcc_2\cup\mcc_3$ and
  $2^{\epsilon(G)-1}> \nu(G)!$ are $Q$-reconstructible as well.
\end{rem}

We end the section with a few open problems.

\begin{problem}
  \label{prob-lovasz} Prove that all graphs $G$ such that
  $\epsilon(G) > \binom{\nu(G)}{2}/2$ are $Q$-reconstructible.
\end{problem}

\begin{problem}
  \label{prob-muller} Prove that if the edge reconstruction conjecture
  is false, then there are infinitely many graphs $G$ such that
  $2^{\epsilon(G)-1}> \nu(G)!$ (preferably with
  $\epsilon(G)=\binom{\nu(G)}{2}/2$) that are not $Q$-reconstructible.
\end{problem}

\begin{problem}
  \label{prob-nw} Counter examples to the edge reconstruction
  conjecture, if they exist, are characterised by a lemma of
  \citet{nw1978}; see also \citet{bondy1991}. Is there a
  characterisation, analogous to the lemma of Nash-Williams, of graphs
  that are not $Q$-reconstructible ?
\end{problem}

\begin{problem}
  Let $\lab_E$ be the class of graphs that are not edge reconstructible,
  and let $\lab_Q$ be the class of graphs that are not
  $Q$-reconstructible. %
  \nomenclature[4g1e]{$\lab_E$}{class of graphs that are not edge
    reconstructible}%
  \nomenclature[4g1q]{$\lab_Q$}{class of graphs that are not
    $Q$-reconstructible}%
  We have shown that $\mcc_1\cup\mcc_2\cup\mcc_3 \subseteq \lab_Q$,
  where equality holds if the edge reconstruction conjecture is true.
  If the reconstruction conjecture is false, then we only know that
  $\lab_E \subseteq \lab_Q$ (ignoring isolated vertices in graphs in
  $\lab_E$), but we do not know if there are graphs in
  $\lab_Q \setminus (\lab_E \cup \mcc_1\cup\mcc_2\cup\mcc_3)$. We showed
  in \cite{thatte2005} that if Ulam's conjecture is false, and if $G$
  and $H$ are non-isomorphic graphs with the same deck, then
  $\isp{2G} = \isp{2H}$. Is there an analogous result for the edge
  reconstruction problem?  %
\end{problem}

\section{Homomorphism cancellation}
\label{sec-hom}
Let $G, H \in \lab$. A \emph{homo\-morphism} from $G$ to $H$ is a map
$f:V(G)\to V(H)$ such that if $\{x,y\}$ is an edge in $G$ then
$\{f(x),f(y)\}$ is an edge in $H$. A one-to-one homomorphism is called a
{\em monomorphism}. Let $\hom(G,H)$ denote the number of homomorphisms
from $G$ to $H$, and let $\mon(G,H)$ denote the number of monomorphisms
from $G$ to $H$. %
\nomenclature[4hom]{$\hom(G,H)$}{number of homomorphisms from $G$ to
  $H$} %
\nomenclature[4mon]{$\mon(G,H)$}{number of monomorphisms from $G$ to
  $H$} %
Both these parameters are well-defined even when one or both of $G$ and
$H$ is unlabelled, since $\hom(G,H) = \hom(G^\prime,H^\prime)$ and
$\mon(G,H) = \mon(G^\prime,H^\prime)$ whenever $G \cong G^\prime$ and
$H \cong H^\prime$. Given an unlabelled graph $G$, we denote by
$\rep{G}$ a representative labelled graph in $G$. %
\citet{lovasz1971} proved the following result.

\begin{thm}[\citet{lovasz1971}; see also Problem 20, Chapter 13 in
  \citet{lovasz1993}]
  \label{thm-lov}
  Let $G_1, G_2 \in \lab$.
  \begin{enumerate}
  \item If $\hom(G_1,H) = \hom(G_2,H)$ for all $H \in \lab$, then
    $G_1 \cong G_2$.
  \item If $\hom(H,G_1) = \hom(H,G_2)$ for all $H \in \lab$, then
    $G_1 \cong G_2$.
\end{enumerate}
\end{thm}

We propose the following conjecture, which in a sense generalises the
idea of homomorphism cancellation in Theorem~\ref{thm-lov}.

\begin{con}
  \label{conj-lov}
  Let $\pi:(\unlab)\bij (\unlab)$ be a bijection such that
  $\hom(G,H) = \hom(\pi(G),\pi(H))$ for all $G, H \in \unlab$. Then
  $G = \pi(G)$ for all $G \in \unlab$.
\end{con}

We show in Proposition~\ref{prop-hom2} that Conjecture~\ref{conj-lov} is
weaker than the edge reconstruction conjecture.

\begin{lem}
  \label{lem-hom1}
  Let $\pi:(\unlab)\to (\unlab)$ be a bijection such that
  $\hom(G,H) = \hom(\pi(G),\pi(H))$ for all $G, H \in \unlab$.  Then
  $\nu(G) = \nu(\pi(G))$ and $\epsilon(G) = \epsilon(\pi(G))$, and
  $\pi(G) = G$ for all $G$ such that $\epsilon(G) \leq 3$.
\end{lem}

\begin{proof}
  \begin{pclaims}

    \pclaim\label{c-phi} $\pi(\Phi) = \Phi$ (where $\Phi$ denotes the
    null graph).
    
    \ppf Let $\pi(G) = \Phi $ and $\pi(\Phi) = H$ for some graphs $G$
    and $H$. Therefore,
    $\hom(G,\Phi) = \hom(\pi(G),\pi(\Phi)) = \hom(\Phi,H)$. We have
    $\hom(\Phi,H) = 1$, and $\hom(G,\Phi)=1$ if and only if
    $\nu(G) = 0$, i.e., $G = \Phi$. Therefore, $\pi(\Phi) = \Phi$. \pqed

    \pclaim\label{c-k1} $\pi(K_1) = K_1$.

    \ppf Let $\pi(G) = K_1 $ and $\pi(K_1) = H$ for some graphs $G$ and
    $H$. Since, $\pi(\Phi) = \Phi$, the graphs $G$ and $H$ are
    non-null. Now,
    $\hom(G,K_1)=\hom(\pi(G),\pi(K_1)) = \hom(K_1 ,H) = \nu(H) \geq 1$.
    We have $\hom(G,K_1) = 1$ if $\epsilon(G) = 0$, and
    $\hom(G,K_1) = 0$ otherwise.  Therefore, $\nu(H) = 1$ and
    $\epsilon(G) = 0$. Moreover, since $H$ is simple, we have
    $\pi(K_1) = H = K_1$.  \pqed

    \pclaim\label{c-nu} For all $G$, we have $\nu(G) = \nu(\pi(G))$.

    \ppf We have
    $\nu(G) = \hom(K_1,G) = \hom(\pi(K_1),\pi(G)) = \hom(K_1,\pi(G)) =
    \nu(\pi(G))$. \pqed

    \pclaim\label{c-k2} $\pi(K_2) = K_2$ and $\pi(2K_1) = 2K_1$.

    \ppf If the claim is not true, then by Claim~\ref{c-nu} we must have
    $\pi(K_2) = 2K_1$ and $\pi(2K_1) = K_2$, but that is not possible
    since $\hom(K_2,K_2) \neq \hom(2K_1,2K_1)$. \pqed

    \pclaim\label{c-eps} For all $G$, we have
    $\epsilon(G) = \epsilon(\pi(G))$.

    \ppf We have
    $2\epsilon(G) = \hom(K_2,G) = \hom(K_2,\pi(G)) =
    2\epsilon(\pi(G))$. \pqed

    \pclaim $\pi(K_3) = K_3$ and $\pi(K_{1,2}) = K_{1,2}$.

    \ppf The claim follows from Claims~\ref{c-nu} and~\ref{c-eps}, and
    that every graph $G$ on at most 3 vertices is determined by the pair
    $(\nu(G), \epsilon(G))$. \pqed

    \pclaim If $\epsilon(G) \leq 3$ then $\pi(G) = G$.

    \ppf If $G$ contains a triangle, then $\pi(G) = G$, which follows
    from Claim~\ref{c-nu} and that $\hom(K_3,G)=\hom(K_3,\pi(G))$.  If
    $G$ does not contain a triangle and has at most 3 edges, then $G$ is
    uniquely determined by the triple
    $\nu(G), \epsilon(G), \hom(K_{1,2},G)$.  \qedhere
  \end{pclaims} 
\end{proof}

\begin{prop}
  \label{prop-hom2}
  The edge reconstruction conjecture implies Conjecture~\ref{conj-lov}.
\end{prop}

\begin{proof}
  Let $\pi:(\unlab)\to (\unlab)$ be a bijection such that
  $\hom(G,H) = \hom(\pi(G),\pi(H))$ for all $G, H \in \unlab$. We assume
  the edge reconstruction conjecture to be true, and prove by induction
  on the number of edges that $\pi(G) = G$ for all $G$. In
  Lemma~\ref{lem-hom1}, we proved that $\pi(G) = G$ for all $G$ such
  that $\epsilon(G)\leq 3$.  Suppose that $\pi(G) = G$ for all graphs
  $G$ such that $3 \leq \epsilon(G) \leq m$. Let $G \in \unlab$ be a
  graph with $m+1$ edges.

  For an unlabelled graph $H$ and an equivalence relation $\Theta$ on
  $V(\rep{H})$, let $\rep{H}\slash\Theta$ denote the graph obtained by
  identifying vertices in each equivalence class of $\Theta $.  Each
  homomorphism from $\rep{H}$ to $\rep{G}$ is a monomorphism from
  $\rep{H}\slash\Theta$ to $\rep{G}$ for some equivalence relation
  $\Theta $ on $V(\rep{H})$.  Therefore,

  \begin{equation}
    \label{eq-hom1}
    \hom(\rep{H},\rep{G}) = \sum_{\Theta} \mon(\rep{H}\slash\Theta,\rep{G}),
  \end{equation}
  where the summation is over all equivalence relations on $V(\rep{H})$.
  In general, for all equivalence relations $\Theta$ on $V(\rep{H})$ we
  have
  \begin{equation}
    \label{eq-hom2}
    \hom(\rep{H}\slash\Theta,\rep{G}) = 
    \sum_{\Theta^\prime \,\mid\, \Theta \leq \Theta^\prime}
    \mon(\rep{H}\slash\Theta^\prime,\rep{G}),
  \end{equation}
  where $\Theta\leq \Theta^\prime$ means $\Theta$ is a refinement of
  $\Theta^\prime$. Following \citet[Chapter 15, Problem
  20]{lovasz1993}, we solve the system of Equations~(\ref{eq-hom2}) for
  $\mon(\rep{H},\rep{G})$ in terms of
  $\hom(\rep{H}\slash\Theta,\rep{G})$, and write
  \begin{equation}
    \label{eq-mon}
    \mon(H,G) = \mon(\rep{H},\rep{G}) = \sum_{\Theta} \alpha_{\rep{H}\slash\Theta} 
    \hom(\rep{H}\slash\Theta,\rep{G}),
  \end{equation}
  where $\alpha_{\rep{H}\slash\Theta}$ are constants (that do not depend
  on $G$). (Another way to look at the solutions of the system of
  equations is via M\"obius inversion.)

  For all $H \in \unlab$ such that $\epsilon(H) \leq m$, and for all
  equivalence relations $\Theta$ on $V(\rep{H})$, we have, by induction
  hypothesis,
  $\hom(\rep{H}\slash\Theta,G) =
  \hom(\rep{H}\slash\Theta,\pi(G))$.
  Hence $\mon(H,G) = \mon(H,\pi(G))$ (by Equation~\ref{eq-mon}). In
  other words, $G$ and $\pi(G)$ have the same edge-deck. Now the edge
  reconstruction conjecture implies that $G = \pi(G)$, completing the
  induction step, and the result.
\end{proof}

The following statement is analogous to Conjecture~\ref{conj-lov}, but
for labelled graphs.
\begin{con}
  \label{conj-lov1}
  Let $\pi:\lab\to \lab$ be a bijection such that
  $\hom(G,H) = \hom(\pi(G),\pi(H))$ for all $G, H \in \lab$. Then
  $G \cong \pi(G)$ for all $G \in \lab$.
\end{con}
It is unclear if Conjecture~\ref{conj-lov} and
Conjecture~\ref{conj-lov1} are equivalent, although it is tempting to
believe that they are. The edge reconstruction conjecture implies
Conjecture~\ref{conj-lov1} as well; the proof of this fact is similar to
the proof of Proposition~\ref{prop-hom2}, and we skip it.

\subsection*{Acknowledgements} 
I would like to thank the following institution for their support: Allan
Wilson Centre for Molecular Ecology and Evolution, New Zealand
(2005-2006), Alfr\'ed R\'enyi Institute of Mathematics, Hungary (2008,
funded by the project Finite Structures (FiSt)), and Universidade de
S\~ao Paulo, Brasil (2011-2013, funded by CNPq, Processo 151782/2010-5
and with additional support from the project MaCLinC).

\printnomenclature[4cm] %

\end{document}